\documentclass[10pt]{amsart}
\usepackage{graphicx}
\usepackage{amsmath,amsthm,amssymb}
\usepackage{mathrsfs}
\usepackage{epstopdf}
\usepackage[T1]{fontenc}
\usepackage[utf8]{inputenc}
\usepackage{enumerate}
\usepackage{physics}
\usepackage{url}
\usepackage{bbm}
\usepackage{hyperref}
\usepackage{array}
\usepackage{todonotes}
\usepackage[sort,numbers]{natbib}
\usepackage{multibib}
\newcites{New}{Additional references}

\usepackage{url}
\makeatletter
\g@addto@macro{\UrlBreaks}{\UrlOrds}
\makeatother

\providecommand{\noopsort}[1]{} 

\newcommand{\tend}[3][]{\xrightarrow[#2\to#3]{#1}}
\newcommand{\Z}{\mathbb{Z}}
\newcommand{\R}{\mathbb{R}}
\newcommand{\N}{\mathbb{N}}
\newcommand{\A}{\mathbb{A}}
\newcommand{\D}{\mathbb{D}}
\newcommand{\U}{\mathbb{U}}
\newcommand{\C}{\mathbb{C}}
\newcommand{\PP}{\mathbb{P}}
\newcommand{\vep}{\varepsilon}
\newcommand{\mob}{\boldsymbol{\mu}}
\newcommand{\bfu}{\boldsymbol{u}}
\newcommand{\bfv}{\boldsymbol{v}}
\newcommand{\lio}{\boldsymbol{\lambda}}
\newcommand{\ot}{\otimes}
\newcommand{\ov}{\overline}

\makeatletter
\def\paragraph{\@startsection{paragraph}{4}%
  \z@\z@{-\fontdimen2\font}%
  {\normalfont\bfseries}}
\makeatother

\newtheorem{Th}{Theorem}[section]
\newtheorem{Prop}[Th]{Proposition}

\theoremstyle{definition}
\newtheorem{Remark}[Th]{Remark}

\newtheorem{Cor}[Th]{Corollary}

\newtheorem{Question}{Question}

\title{Sarnak's Conjecture from the Ergodic Theory Point of View}

\author{Joanna Ku\l{}aga-Przymus \and Mariusz Lema\'{n}czyk}

\begin{document}

\maketitle

\section{Glossary}
\paragraph{Measure-theoretic dynamical system}
Let $(Z,\mathscr{C},\kappa)$ be a standard probability Borel space. Assume $R\colon X\to X$ is invertible (a.e.), bi-measurable  and $\kappa$-preserving. Then $R$ is called an {\em automorphism} or (invertible) {\em measure-preserving transformation} and the quadruple $(R,Z,\mathscr{C},\kappa)$ is a {\em measure-theoretic dynamical system}.
\paragraph{Topological dynamical system}
Let $X$ be a compact metric space and let $T$ be a homeomorphism of $X$. Then $(T,X)$ is called a {\em topological dynamical system}.

\paragraph{Subshift}
For any closed subset $\mathbb{A}\subset \mathbb{U}:=\{z \in \mathbb{C} : |z|\leq 1\}$ (most often $\mathbb{A}$ will be finite), let $S\colon \A^\Z\to \A^\Z$ be the left shift, i.e. for each $x=(x_n)_{n\in\Z}\in \A^\Z$, $Sx=y$, where $y_n=x_{n+1}$ for $n\in\Z$. Each closed and $S$-invariant subset $X\subset \A^\Z$ is called a {\em subshift} and the corresponding dynamical system is $(S,X)$. For each $u\in \A^\Z$, we define a subshift $X_u\subset \A^{\Z}$ as the orbit closure of $u$ under $S$. Similarly, when $u\in \A^\N$ {($\N=\{1,2,\ldots\}$)}, we can extend $u$ symmetrically, by setting $u(-n):=u(n)$ {(and $u(0)=0$)} for each $n\geq 1$ and define again the corresponding subshift $X_u$.
Finally, we will denote by $F$ the continuous function on $X\subset\mathbb{A}^\Z$, defined by $F(x)=x_0$.

\paragraph{Invariant measure}
Given a topological dynamical system $(T,X)$, the set of probability Borel $T$-invariant measures is denoted by $M(T,X)$. The subset of ergodic measures (which is always non-empty) is denoted by $M^e(T,X)$. Each $\nu\in M(T,X)$ gives rise to a measure-theoretic dynamical system $(T,X,\mathscr{B}(X),\nu)$, where $\mathscr{B}(X)$ stands for the $\sigma$-algebra of Borel subsets of $X$. With the weak-$*$-topology, $M(T,X)$ becomes a compact metrizable space. $(T,X)$ is called {\em uniquely ergodic} if $|M(T,X)|=1$.

\paragraph{Uniquely ergodic model of a measure-theoretic dynamical  system}
Given a measure-theoretic (ergodic) dynamical system $(R,Z,\mathscr{C},\kappa)$ and a uniquely ergodic topological system $(T,X)$ with $M(T,X)=\{\nu\}$ ($\nu$ is necessarily ergodic), one says that $(T,X)$ is a {\em uniquely ergodic model} of the automorphism $R$ if the measure-theoretic systems $(R,Z,\mathscr{C},\kappa)$ and $(T,X,\mathscr{B}(X),\nu)$ are measure-theoretically isomorphic. Each ergodic automorphism $R$ has a uniquely ergodic model.
\paragraph{Generic point}
Assume that $(T,X)$ is a topological dynamical system.
We say that $x\in X$ is a {\em generic point} for a Borel measure $\nu$ on $X$, whenever the ergodic theorem holds for $T$ at $x$ for any continuous function $f\in C(X)$, i.e.\ $\frac{1}{N}\sum_{n\leq N}f(T^nx)\to \int f\, d\nu$. In other words, the empirical measures $\frac{1}{N}\sum_{n\leq N}\delta_{T^nx}$ converge to $\nu$ in the weak-$*$-topology. If this convergence holds along a subsequence, we say that $x$ is {\em quasi-generic} for $\nu$. We denote by $V(x)$ the set of all $\nu\in M(T,X)$ for which $x$ is quasi-generic ($\emptyset\neq V(x)\subset M(T,X)$ by compactness).

We say that $x\in X$ is {\em logarithmically generic} for $\nu$, whenever $\frac{1}{{L_N}}\sum_{n\leq N}\frac{1}{n}\delta_{T^{n}x}$, {with $L_N=\sum_{n\leq N}1/n$}, converges to $\nu$. {No harm arises if in what follows we replace $L_N$ by $\log N$.} Finally, we say that $x\in X$ is {\em logarithmically quasi-generic} for $\nu$ whenever this convergence holds along a subsequence and we denote the set of all measures for which $x$ is logarithmically quasi-generic by $V^{\log}(x)$, ($\emptyset\neq V^{\log}(x)\subset M(T,X)$).

\paragraph{{Quasi-genericity versus logarithmic quasi-genericity}} {Obviously, if $|V(x)|=1$ then $V(x)=V^{\log}(x)$, but in general the sets $V(x)$ and $V^{\log}(x)$ can be even disjoint. However,} as shown in~\cite{MR3821718},
\begin{equation}\label{gkl1}
V^{\log}(x)\cap M^e(T,X)\subset V(x).
\end{equation}
Moreover, using an idea from Tao~\cite{Ta5}, it has been proved in \cite{Gomilko:ab} that
\begin{multline}\label{glr1}
\mbox{If $V^{\log}(x)=\{\nu\}$ and $\nu$ is ergodic, then $\lim_{k\to\infty}\frac1{N_k}\sum_{n\leq N_k}\delta_{T^nx}=\nu$},\\
\mbox{for a subset $\{N_k:k\geq1\}$ whose logarithmic density is~1.}
\end{multline}

\paragraph{Entropy}
There are two basic notions of {\em entropy}: topological and measure-theoretic. We skip the definitions and refer the reader, e.g., to \cite{MR2809170}. For a topological dynamical system $(T,X)$ its topological entropy will be denoted by $h_{top}(T,X)$ and for a measure-theoretic dynamical system $(T,X,\mathscr{B},\nu)$ the corresponding measure-theoretic entropy will be denoted by $h(T,X,\mathscr{B},\nu)$. The basic connection between them is the variational principle: $$h_{top}(T,X)=\sup_{\nu\in M(T,X)}h(T,X,\mathscr{B}(X),\nu)=\sup_{\nu\in M^e(T,X)}h(T,X,\mathscr{B}(X),\nu).$$

\paragraph{Completely deterministic point}
We say that point $x\in X$ is {\em completely deterministic}~\cite{We9} (see also~\cite{Kam}) if for any $\nu\in V(x)$, we have $h(T,X,\mathscr{B}(X),\nu)=0$. By the variational principle, $h_{top}(T,X)=0$ if and only if all points of $X$ are completely deterministic.

\paragraph{Furstenberg system}
Let $u\in \A^\Z$. For each $\nu\in V(u)$, the system $(T,X_u,\mathscr{B}(X_u),\nu)$ is called a {\em Furstenberg system of $u$}.
For each $\nu\in V^{\log}(u)$, the system $(T,X_u,\mathscr{B}(X_u),\nu)$ is called a {\em logarithmic Furstenberg system of $u$}.

For each $u\in \A^\Z$, one can consider $|u|\in [0,1]^\Z$. Then $(S,X_{|u|})$ is a topological factor of $(S,X_u)$ (the map $\pi\colon X_{u}\to X_{|u|}$ given by $\pi(x)=|x|$, understood coordinatewise, is equivariant with $S$). For $\nu\in V(u)$, we have $\pi_\ast(\nu)\in V(|u|)$, where $\pi_\ast(\nu)$ stands for the image of $\nu$ via $\pi$. Moreover, the Furstenberg system $(S,X_u,\mathscr{B}(X_u),\nu)$ is an extension of $(S,X_{|u|},\mathscr{B}(X_{|u|}),\pi_\ast(\nu))$. In particular, if $V(|u|)$ is a singleton ($|u|$ is a generic point), then all Furstenberg systems of $u$ have $(S,X_{|u|},\mathscr{B}(X_{|u|}),\nu)$ (where $\nu$ is the unique member of $V(|u|)$) as their factor.

\paragraph{Arithmetic function} A sequence of complex numbers is usually denoted by $u=(u_n)$. But if such a sequence is, in some sense, important from number theory point of view, one speaks about an {\em arithmetic function} and rather writes 
$\bfu\colon \N\to \C$, $\bfu=(\bfu(n))$.  

An arithmetic function $\bfu$ is said to be {\em multiplicative}, whenever $\bfu(1)=1$ and $\bfu(m\cdot n)=\bfu(m)\cdot \bfu(n)$ for any choice of coprime $m,n\in \N$. The prominent examples of multiplicative functions are the M\"obius function $\mob$ and the Liouville function $\lio$. The M\"obius function $\mob\colon\N\to\{-1,0,1\}$ is defined by $\mob(p_1\ldots p_k)=(-1)^k$ for different prime numbers $p_1,\ldots,p_k$ (in what follows, the set of primes is denoted by $\PP$), $\mob(1)=1$ and $\mob(n)=0$ for all non-square-free numbers. The Liouville function $\lio\colon\N\to\{-1,1\}$ is given by $\lio(n)=(-1)^{i_1+\ldots+i_k}$ for $n=p_1^{i_1}\ldots p_k^{i_k}$ with $p_1,\ldots,p_k\in\PP$ and $i_1,\ldots,i_k\in\N$. Clearly $\mob=\lio\cdot\mob^2$, where $\mob^2$ is nothing but the characteristic function of the set $\mathscr{S}$ of square-free numbers. In fact, $\lio$ is completely multiplicative, i.e.\ $\lio(m\cdot n)=\lio(m)\cdot \lio(n)$ for any choice of $m,n\in \N$. We extend both, $\mob$ and $\lio$, to negative coordinates symmetrically.

\paragraph{Aperiodicity}
We say that $\bfu\colon \N\to \C$ is {\em aperiodic} whenever $\bfu$ has a mean, equal to zero, along each arithmetic progression: $\lim_{N\to \infty}\frac{1}{N}\sum_{n\leq N}\bfu(an+b)=0$. Many classical multiplicative functions are aperiodic, including $\mob$ and $\lio$.

A distance between $\bfu,\bfv\colon \N\to\U$ is defined as
$$\D(\bfu,\bfv;N):=\left(\sum_{p\in\PP,p\leq N}\frac{1-\Re(\bfu(p)\overline{\bfv(p)})}{p} \right)^{1/2}.$$
We say that $\bfu\colon \N\to \U$ is {\em strongly aperiodic}~\cite{MR3435814}, whenever $M(\bfu\cdot \chi; N):=\min_{|t|\leq N}\D(\bfu\cdot \chi,n^{it};N)^2\to \infty$ as $N\to \infty$ for every Dirichlet character $\chi$ (i.e.\ for every periodic, completely multiplicative function). Strong aperiodicity implies aperiodicity. The converse is not in general true~(see Theorem B.1 in~\cite{MR3435814}), but it is true for (bounded) real valued multiplicative functions (see Appendix C in~\cite{MR3435814}). In particular, $\mob$ and $\lio$ are strongly aperiodic.

\paragraph{Orthogonality of sequences}
Suppose that one of sequences $(u_n),(v_n)\subset\C$ is of zero {mean}. We say that $(u_n),(v_n)$ are {\em orthogonal}, whenever $\lim_{N\to\infty}\frac{1}{N}\sum_{n\leq N}u_n \overline{v}_n=0$. We say that sequences $(u_n),(v_n)\subset\C$ are {\em orthogonal on short intervals}, whenever $\lim_{K\to\infty}\frac{1}{b_K}\sum_{k\leq K}\left| \sum_{b_k\leq n<b_{k+1}}u_n\overline{v}_n\right|= 0$
for any sequence $(b_k)$ of natural numbers with $b_{k+1}-b_k\to \infty$. Clearly, orthogonality on short intervals implies orthogonality. We can consider $u_n$ being an element of a Banach space, while in the dynamical context one often takes $u_n=f(T^nx_k)$, whenever  $n\in [b_k, b_{k+1})$ (with ${x_k}\in X$ and $f\in C(X)$). Moreover, often $(v_n)$ will be in fact a multiplicative function (and then we rather write $(\bfv(n))$).

\paragraph{M\"obius orthogonality}
We say that a (topological) dynamical system $(T,X)$ is \emph{M\"obius orthogonal} if
\begin{equation}\label{defmo}
\lim_{N\to\infty}\frac{1}{N}\sum_{n\leq N}f(T^nx)\mob(n)=0
\end{equation}
for each $f\in C(X)$ and $x\in X$.

\paragraph{Joinings of measure-theoretic dynamical systems} Assume that $(R_i,Z_i,\mathscr{C}_i,\kappa_i)$ is a measure-theoretic dynamical system, $i=1,2$. Each $R_1\times R_2$-invariant measure $\rho$ on $\mathscr{C}_1\otimes\mathscr{C}_2$ projecting on $\kappa_1$ and $\kappa_2$, respectively, is called a {\em joining} of the automorphisms $R_1$ and $R_2$. The set of joinings between $R_1$ and $R_2$ is denoted by $J(R_1,R_2)$ and each $\rho\in J(R_1,R_2)$ yields a (new) measure-theoretic dynamical system $(R_1\times R_2, Z_1\times Z_2,\mathscr{C}_1\otimes\mathscr{C}_2,\rho)$. When $R_1,R_2$ are both ergodic then the set $J^e(R_1,R_2)$ of ergodic joinings between $R_1$ and $R_2$ is non-empty. If $J(R_1,R_2)=\{\kappa_1\otimes\kappa_2\}$ then $R_1$ and $R_2$ are called {\em disjoint} (in the sense of Furstenberg).

Let $(T,X)$ be a topological dynamical system. Fix $x\in X$ and let $\bfu$ be an arithmetic function bounded by $1$, i.e.\ $\bfu\colon \N\to \U$. Each accumulation point $\rho$ of $\frac{1}{N}\sum_{n\leq N} \delta_{(T^nx,S^n\bfu)}$ is a $(T\times S)$-invariant measure on $X\times \U^\Z$. Obviously, it is a {\em joining} of its projections onto both coordinates. Hence, $\rho$ is a joining of $(T,X,\mathscr{B}(X),\rho|_X)$ with a Furstenberg system of $\bfu$.

\paragraph{Nilrotation}
Let $G$ be a connected, simply connected nilpotent Lie group and $\Gamma\subset G$ a lattice (a discrete, cocompact subgroup).
For any $g_0\in G$ we define $T_{g_0}(g\Gamma):=g_0g\Gamma$. Then the topological system $(T_{g_0},G/\Gamma)$ is called a {\em nilrotation}.

\section{Definition of the Subject}
Sarnak in 2010 formulated a now celebrated conjecture on the M\"obius orthogonality. It states that each topological dynamical system $(T,X)$ of zero entropy is M\"obius orthogonal, i.e.\ whenever $h_{top}(T,X)=0$,  
\begin{equation}\label{e1}
\lim_{N\to\infty}\frac{1}{N}\sum_{n\leq N}f(T^nx)\mob(n)=0\text{ for each } f\in C(X) \text{ and }x\in X.
\end{equation}
It has become a bridge between analytic number theory and dynamics. In particular, it keeps stimulating a quick development of disjointness theory in dynamics originated by Furstenberg in 1967, with potential applications, for example in number theory. In this article we will focus on ergodic-theoretic aspects of recent progress in the area. {Unless we need this for a historical reason, we rather avoid describing numerous particular classes of topological systems in which M\"obius disjointness
has been shown (an extended bibliography is provided at the end of the article). We concentrate on main ideas, and describe some general results.}

\section{Introduction}
The article is organized as follows. In Section~\ref{se:chowla} and Section~\ref{se:sarnak} respectively, we discuss Chowla and Sarnak's conjecture from the ergodic-theoretic viewpoint, including dynamical interpretations of purely number-theoretic statements and the main strategies used to attack Sarnak's conjecture. Section~\ref{se:arithm} is a survey of results on Sarnak's conjecture, arranged with respect to the properties of the M\"obius function that come into play in the proof. Finally, in Section~\ref{se:future}, we state some open problems.

\section{Chowla Conjecture}\label{se:chowla}

\subsection{\eqref{chsfor}}
The Chowla conjecture deals with  higher order correlations of the M\"obius function and asserts that
\begin{equation}\label{chsfor}\tag{{\bf C}}
\lim_{N\to \infty}\frac1N\sum_{n\leq N}\mob^{j_0}(n+k_0)\mob^{j_1}(n+k_1)\ldots\mob^{j_r}(n+k_r)=0
\end{equation}
whenever $0\leq k_0<\ldots<k_r$, $j_s\in\{1,2\}$ not all equal to 2, $r\geq0$.\footnote{In \cite{Ch}, the Chowla conjecture it is formulated for the Liouville function. We follow \cite{Sa}. For a discussion on an equivalence of the Chowla conjecture with $\mob$ and $\lio$, see~\cite{MR3821722}. } The above statement can be translated into dynamical language. To see this, consider first $|\mob|=\mob^2$, i.e.\ the characteristic function of the set $\mathscr{S}$ of square-free numbers. As a member of $\{0,1\}^\Z$, it is well-known to be a generic point for {an ergodic measure, namely for} the so-called Mirsky measure $\nu_{\mathscr{S}}$ {(also denoted as $\nu_{\mob^2}$)}, given by the ordinary average frequencies of blocks, see e.g.\ \cite{Ab-Ku-Le-Ru}. \eqref{chsfor} is equivalent to the following:
\begin{equation}
    \begin{minipage}{.8\textwidth}
        \centering
	the point $\mob\in\{-1,0,1\}^{\Z}$ is generic for the relatively independent extension $\widehat{\nu}_{\mathscr{S}}$ of $\nu_{\mathscr{S}}$, via the natural map $\{-1,0,1\}^{\Z}\ni(x_n)_{n\geq1}{\mapsto}(x^2_n)_{n\geq1}\in\{0,1\}^{\Z}$.
    \end{minipage}
\end{equation}
The measure $\widehat{\nu}_{\mathscr{S}}$ is given by the following condition: for each block $C$ over the alphabet $\{-1,0,1\}$, we have $\widehat{\nu}_{\mathscr{S}}(C)=2^{-k}\nu_{\mathscr{S}}(C^2)$, where $C^2$ is obtained from $C$ by taking the square (or, equivalently, {absolute value}) at each coordinate and $k$ is the cardinality of the support of $C$. To see this, we use the following:
\begin{Remark}
The {span of the} family of continuous functions $\{F^{j_0}\circ S^{k_0} \cdot F^{j_1}\circ S^{k_1}\cdot \ldots\cdot F^{j_\ell}\circ S^{k_\ell} : j_i \geq 0, k_i\in \Z\}$ forms an algebra that distinguishes points. It follows directly from the Stone-Weierstrass theorem that the values of  integrals $\int F^{k_1}\circ S^{r_1} \cdot F^{k_2}\circ S^{r_2}\cdot \ldots\cdot F^{k_\ell}\circ S^{r_\ell} \, d\kappa$ determine measure $\kappa$.
\end{Remark}
Now, it suffices to look at measures $\kappa\in V(\mob)$ and compare the value of the integrals $\int F^{k_1}\circ S^{r_1} \cdot F^{k_2}\circ S^{r_2}\cdot \ldots\cdot F^{k_\ell}\circ S^{r_\ell} \, d\kappa$ with the corresponding values of  $\int F^{k_1}\circ S^{r_1} \cdot F^{k_2}\circ S^{r_2}\cdot \ldots\cdot F^{k_\ell}\circ S^{r_\ell} \, d\widehat{\nu}_{\mathscr{S}}$. We will also use the fact that
\begin{equation}\label{ergmu}
\widehat\nu_{\mathscr{S}}\in M^e(S,X_{\mob}).
\end{equation}

The simplest instance of Chowla conjecture, i.e.\ $\lim_{N\to\infty}\frac{1}{N}\sum_{n\leq N}\mob(n)=0$, is known to hold and it is equivalent to the Prime Number Theorem, as shown by Landau. Moreover, if exactly one of the exponents $j_i$ is odd, then~\eqref{chsfor} also holds (in~\cite{MR3810678} a quantitative version of this fact has been proved). We also have the following conditional result of Frantzikinakis on~\eqref{chsfor}:
\begin{Th}[\cite{MR3742396}]\label{t:fr1}
Assume that $V(\lio)=\{\kappa\}$ and $\kappa$ is ergodic. Then~\eqref{chsfor} holds for $\lio$.\footnote{\label{f:nikos} That is, $\kappa$ equals the $\frac12$-Bernoulli measure on $\{-1,1\}^{\Z}$. Cf.\ Theorem~\ref{t:fr1a}.}
\end{Th}

\subsection{\eqref{chsfor} vs. \eqref{clog}}
Analogously, one can study the logarithmic Chowla conjecture~\cite{MR3676413}, asserting that
\begin{equation}\tag{{\bf C$_{\log}$}}\label{clog}
\lim_{N\to\infty}\frac1{\log N}\sum_{n\leq N}\frac1n\mob^{j_0}(n+k_0)\mob^{j_1}(n+k_1)\ldots\mob^{j_r}(n+k_r)=0,
\end{equation}
with all parameters as before. Again, this can be translated to the dynamical language, using the notion of a logarithmically generic point.
Thus, using \eqref{ergmu},~\eqref{gkl1} and~\eqref{glr1}, we have the following:
\begin{Cor}[cf.\ \cite{MR3821718,MR3676413}]\label{wnio}
\eqref{clog} implies \eqref{chsfor} along a subsequence of full logarithmic density.
\end{Cor}

\subsection{More on \eqref{clog}}
We begin this section with a conditional result of Frantzikinakis on~\eqref{clog}:\footnote{Cf.\ Footnote~\ref{f:nikos}.}
\begin{Th}[\cite{MR3742396}]\label{t:fr1a}
Assume that $\kappa\in V^{\log}(\lio)$ is an ergodic measure. Then \eqref{clog} holds for $\lio$ along the same subsequence for which $\lio$ is quasi-generic for $\kappa$.
\end{Th}
The above result is already formulated in the language of (logarithmically) quasi-generic points. We pass now to number-theoretic results, with their ergodic-theoretic consequences.
\begin{Th}[\cite{MR3569059}]\label{ch2}
\eqref{clog} holds for $r=1$ (for $\lio$ in place of $\mob$), i.e.\ for each $0\neq h\in\Z$, we have
$$
\lim_{N\to\infty}\frac1{\log N}\sum_{n\leq N}\frac{\lio(n)\lio(n+h)}n=0.
$$
\end{Th}
In fact, Tao in~\cite{MR3569059} proves a stronger result, being an instance of logarithmically averaged Elliott conjecture. One of the consequences is that the analogue of Theorem~\ref{ch2} holds for $\mob$:
\begin{Cor}[\cite{MR3569059}]\label{mob2}
\eqref{clog} holds for $r=1$, i.e.\ for each $0\neq h\in\Z$, we have
$$
\lim_{N\to\infty}\frac1{\log N}\sum_{n\leq N}\frac{\mob(n)\mob(n+h)}n=0.
$$
\end{Cor}
This, in turn, has the following interpretation in terms of spectral measures:
\begin{Cor}[\cite{MR3821717}]\label{lebesgue}
For each logarithmic Furstenberg system $(S,X_{\mob},\mathscr{B}(X_{\mob}),\kappa)$ of $\mob$, the spectral measure $\sigma_F$ of $F$ is Lebesgue. The same holds for $\kappa\in V^{\log}(\lio)$.
\end{Cor}

Moreover, we have the following:
\begin{Th}[\cite{MR3992031,MR3938639}]\label{odd corr}
\eqref{clog} holds for odd order correlations, i.e.\ for all even values of $r$.
\end{Th}
Again, we want to interpret this from the dynamical viewpoint:
\begin{Cor}
For each logarithmic Furstenberg system $(S,X_{\mob},\mathscr{B}(X_{\mob}),\kappa)$, the element $Ry:=-y$ preserves measure $\kappa$ and commutes with the shift $S$, hence $R$ belongs to the centralizer of $S$.
\end{Cor}
To see that the above statement is true, it suffices to take $\kappa\in V^{\log}(\mob)$ and check that we also have $\kappa\in V^{\log}(-\mob)$. For this, one uses the family of continuous functions $\{F^{j_0}\circ S^{k_0} \cdot F^{j_1}\circ S^{k_1}\cdot \ldots\cdot F^{j_\ell}\circ S^{k_\ell} : j_i \geq 0, k_i\in \Z\}$. Then, for even correlations the desired equalities are obvious and for odd correlations one applies Theorem~\ref{odd corr}.

Finally, as a consequence {of Corollary~\ref{wnio}}, we have the following:
\begin{Cor}[\cite{MR3821718}]
Suppose that \eqref{clog} holds. Then $\hat{\nu}_{\mob^2}\in V(\mob)$.
\end{Cor}

\subsection{Averaged~\eqref{chsfor}}
Also, an averaged version of Chowla conjecture is present in the literature. As we will see later, this form of averaging will play a special role from the point of view of Sarnak's program, under the name of convergence on short intervals. Matom\"aki, Radziwi\l\l \ and Tao showed the following result on the order two correlations of $\mob$:
\begin{Th}[\cite{MR3435814}]\label{mrt}
We have
$$
\lim_{\substack{M,H\to \infty\\ \text{ with }H=o(M)}}\frac{1}{HM}\sum_{h\leq H}\left|\sum_{m\leq M}\mob(m)\mob(m+h)\right|=0.
$$
\end{Th}
As a consequence of the above, we have (cf.\ Corollary~\ref{lebesgue}):
\begin{Cor}[\cite{MR3821717}]\label{con1}
For each Furstenberg system $(S,X_{\mob},\mathscr{B}(X_{\mob}),\kappa)$,
the spectral measure $\sigma_F$ of $F$  is continuous.
\end{Cor}

\subsection{\eqref{chsfor} vs. other conjectures}\label{ell}
Chowla conjecture is thought of as a multiplicative analogue of the twin primes problem. Indeed, Twin Primes Conjecture in its quantitative form expects that $\sum_{n\leq N}\boldsymbol\Lambda(n)\boldsymbol\Lambda(n+2)=(2\prod_{\PP\ni p\geq 3}(1-\frac1{(p-1)^2}))N+{\rm o}(N)$, where $\boldsymbol\Lambda$ is the von Mangoldt function.\footnote{We have $\boldsymbol\Lambda(p^k)=\log p$ for each prime $p$ and $k\geq1$, and $\boldsymbol\Lambda$ vanishes at all other values of $n$. This function is a good approximation of $\mathbbm{1}_{\PP}$ and it is {\bf not} a multiplicative function.}

Moreover, Chowla conjecture is a special instance of Elliott conjecture on correlations of multiplicative functions. Let $\bfu_0,\dots, \bfu_k\colon \N\to \U$ be multiplicative. Elliott conjecture, in a corrected form given in~\cite{MR3435814}, asserts that
$$
\lim_{N\to\infty}\frac{1}{N}\sum_{n\leq N}\bfu_0(n+k_0)\bfu_1(n+k_1)\dots \bfu_k(n+k_r)=0,
$$
if for some $0\leq j\leq k$, $\bfu_j$ is strongly aperiodic. This conjecture was stated first in~\cite{MR1222182,MR1292619} and in its original version turned out to be false, see \cite{MR3435814} for details.
Also, a logarithmically averaged version of Elliott conjecture appears in the literature. For the details, see~\cite{MR4039498,MR3992031} and references therein.

\section{Sarnak's Conjecture}\label{se:sarnak}
\subsection{\eqref{sar}}
Sarnak's Conjecture from 2010 states that
\begin{equation}\tag{{\bf S}}\label{sar}
    \begin{minipage}{.9\textwidth}
        \centering
	each (topological) zero entropy system $(T,X)$ is M\"obius orthogonal, i.e.\ satisfies~\eqref{defmo} for all $f\in C(X)$ and $x\in X$.
    \end{minipage}
\end{equation}
As zero entropy expresses the fact that the system is ``deterministic'' (or of ``low complexity''), Sarnak's conjecture captures our expectation that prime numbers behave globally as a random sequence, or, more precisely, that they cannot be predicted by a low-complexity object.
One can relax the entropy assumptions on~$T$ in Sarnak's conjecture in the following way:
\begin{Th}[\cite{Ab-Ku-Le-Ru}]
\eqref{sar} is equivalent to M\"obius orthogonality for each topological dynamical system $(T,X)$, each continuous function $f\in C(X)$ and each completely deterministic point $x\in X$ (i.e.\ \eqref{defmo} holds at $x$ for all $f\in C(X)$).
\end{Th}

\begin{Remark}
The difficulty in Sarnak's conjecture comes from the requirement ``for all $x\in X$''. An a.e.\ version of M\"obius orthogonality is true for {\bf all} dynamical systems. The proof makes use of Davenport's estimate~\eqref{daven} (below),  see~\cite{Sa,Ab-Ku-Le-Ru}.
\end{Remark}

Even though Sarnak's Conjecture is defined in terms of topological dynamics, it can be translated to ergodic-theoretic language. Namely
$$
\frac{1}{N}\sum_{n\leq N}f(T^nx)\mob(n)=\int_{X\times X_{\mob}} f\otimes F\, d\left(\frac{1}{N}\sum_{n\leq N}\delta_{(T^nx,S^n\mob)}\right).
$$
Thus, we need to study the properties of joinings given by the limit points of $\frac1N\sum_{n\leq N}\delta_{(T^nx,S^n\mob)}$.

The simplest case where \eqref{sar} is known to hold is the one-point dynamical system: $\lim_{N\to \infty}\frac{1}{N}\sum_{n\leq N}\mob(n)=0$
is equivalent to the Prime Number Theorem. \eqref{sar} for rotations on finite groups is equivalent to the Dirichlet's Prime Number Theorem. For irrational rotations,~\eqref{sar} follows from an old (quantitative) result of Davenport~\cite{Da}: for an arbitrary $A>0$, 
\begin{equation}\label{daven}
\max_{t \in \mathbb{T}}\left|\displaystyle\sum_{n \leq N}e^{2\pi int}{\boldsymbol{\mu}}(n)\right|\leq C_A\frac{N}{\log^{A}N}\text{ for some }C_A>0\text{ and all }N\geq 2.
\end{equation}
As we will see later, an important role in the research around Sarnak's conjecture is played by nil-systems. Green and Tao obtained the following quantitative version of \eqref{sar}:
\begin{Th}[\cite{Gr-Ta}] Let $G$ be a simply-connected
nilpotent Lie group with a discrete and cocompact subgroup $\Gamma$. Let $p \colon \Z \to G$ be any its polynomial sequence\footnote{I.e.\ $p(n)=a_1^{p_1(n)}\ldots a_k^{p_k(n)}$, where $p_j\colon\N\to\N$ is a polynomial, $j=1,\ldots,k$.} and $f\colon G/\Gamma\to \R$  a Lipschitz function.
Then $$\left|\sum_{n\leq N} f(p(n)\Gamma)\mob(n)\right|={\rm O}_{f,G,\Gamma,A}\left(\frac N{\log^AN}\right)$$
for all $A > 0$.
\end{Th}
In particular, all nilrotations are M\"obius orthogonal.

\subsection{\eqref{sar} vs. \eqref{chsfor}}
Sarnak's Conjecture was originally mainly motivated by Chowla conjecture; we have the following result:
\begin{Th}\label{ChtS}
\eqref{chsfor} implies \eqref{sar}.
\end{Th}
Theorem~\ref{ChtS} is already stated in~\cite{Sa}. In fact, it is a purely ergodic theory claim: we have already noticed that both conjectures have their ergodic theory reformulation and a joining proof of Theorem~\ref{ChtS} can be found in \cite{Ab-Ku-Le-Ru}. The main idea is the following: suppose that $\frac1{N_k}\sum_{n\leq N_k}\delta_{(T^nx,S^n\mob)}\to \rho$. The projection of this joining onto $X$ is a zero entropy measure $\kappa$, whereas the projection onto $X_{\mob}$ equals $\widehat{\nu}_\mathscr{S}$ by Chowla conjecture. Moreover, $(S,X_{\mob},\widehat{\nu}_\mathscr{S})$ has the property of being relative Kolmogorov with respect to its factor $(S,X_{\mob^2},{\nu}_\mathscr{S})$. On the other hand, the restriction of $\rho$ to $X\times X_{\mob^2}$ is of relative zero entropy over $X_{\mob^2}$. This yields relative disjointness of $(S,X_{\mob},\widehat{\nu}_\mathscr{S})$ and $(T\times S, X\times {X_{\mob^2}},\rho|_{X\times {X_{\mob^2}}})$ over their common factor $(S,X_{\mob^2},\nu_{\mathscr{S}})$. To complete the proof, we use the orthogonality of $F$ to $L^2(X_{\mob^2},\nu_\mathscr{S})$.
\begin{Remark}
It still remains open whether \eqref{sar} implies \eqref{chsfor}, see however Remark~\ref{tauw}.
\end{Remark}

In~\cite{HUANG2019827}, M\"obius orthogonality for low complexity systems is discussed. Following~\cite{Ferenczi_1997}, we say that the measure-complexity of $\mu\in M(T,X)$ is weaker than $a=(a_n)_{n\geq 1}$ if
$$
\liminf_{n\to\infty}\frac{\min\{m\geq1: \mu(\bigcup_{j=1}^nB_{d_n}(x_i,\vep))>1-\vep\text{ for some } x_1,\ldots,x_m\in X\}}{a_n}=0
$$
for each $\vep>0$ (here $d_n(y,z)=\frac1n\sum_{j=1}^nd(T^jy,T^jz)$).
\begin{Th}[\cite{HUANG2019827}]
Suppose that \eqref{chsfor} holds for correlations of order 2 (i.e.\ for $r=1$). Then $(T,X)$ is M\"obius {orthogonal} whenever all invariant measures for $(T,X)$ are of complexity weaker than $n$.
\end{Th}
To obtain a non-conditional result, Huang, Wang and Ye use a difficult estimate of Matom\"aki, Radziwiłł and Tao (namely, ``Truncated Elliott on the average'', applied to $\mob$) from~\cite{MR3435814}. The cost to be paid is a further strengthening of the assumptions on the complexity of $(T,X)$.
\begin{Th}[\cite{HUANG2019827}]
Suppose that  all invariant measures of $(T,X)$ are of sub-polynomial complexity, i.e.\ their complexity is weaker than $(n^\tau)_{n\geq 1}$ for each $\tau>0$. Then $(T,X)$ is M\"obius orthogonal.
\end{Th}
See~\cite{Huang:aa} for the most recent application of this result.

Finally, let us point out a consequence of the result on correlations of $\mob$ of order~2. Directly from Corollary~\ref{con1}, we have:
\begin{Cor}\label{AA}
All topological dynamical systems whose all invariant measures yield systems with discrete spectrum are M\"obius orthogonal.\footnote{In the uniquely ergodic case, an earlier and independent proof of this fact was given by Huang, Wang and Zhang~\cite{MR3959363} (for the totally uniquely ergodic case, see~\cite{Ab-Le-Ru2}). The result also follows from~\cite{HUANG2019827}.}
\end{Cor}

\subsection{Strong MOMO property}\label{smomo}
Given an arithmetic function $\bfu$, following~\cite{MR3874857}, we say that $(X,T)$ satisfies the \emph{strong $\bfu$-OMO property} if, for any increasing sequence of integers
$0=b_0<b_1<b_2<\cdots$ with $b_{k+1}-b_k\to\infty$, for any sequence $(x_k)$ of points in $X$, and any $f\in C(X)$, we have
  \begin{equation}
    \label{eq:defMOMOSI}
    \frac{1}{b_{K}} \sum_{k< K} \left|\sum_{b_k\le n<b_{k+1}} f(T^{n-b_k}x_k) \bfu(n)\right| \tend{K}{\infty} 0.
  \end{equation}
If $\bfu=\mob$ we speak about the strong MOMO\footnote{The acronym MOMO stands for M\"obius Orthogonality of Moving Orbits.}  property.
Strong MOMO property was introduced in~\cite{MR3874857} to deal with M\"obius orthogonality of uniquely ergodic models of a given measure-theoretic dynamical system. Moreover, we have:
\begin{Th}[\cite{MR3874857}] \label{th:3wki}
The following conditions are equivalent:
\begin{enumerate}[(i)]
\item All zero entropy systems are M\"obius orthogonal, i.e.\ Sarnak's conjecture holds.
\item For each zero entropy system $(T,X)$, we have $\lim_{N\to\infty}\frac1N\sum_{n\leq N}f(T^nx)\mob(n)=0$ when $N\to\infty$, for each $f\in C(X)$, uniformly in $x\in X$, i.e.\
uniform Sarnak's conjecture holds.
\item
All zero entropy systems enjoy the strong MOMO property.
\end{enumerate}
\end{Th}

By taking $f=1$, we obtain that strong $\bfu$-OMO implies the following:
\begin{equation}\label{eq:Mobius-like}
 \frac{1}{b_{K}} \sum_{k< K} \left|\sum_{b_k\le n<b_{k+1}} \bfu(n)\right| \tend{K}{\infty} 0
 \end{equation}
 for every sequence
$0=b_0<b_1<b_2<\cdots$ with $b_{k+1}-b_k\to\infty$.
In particular, $\frac1N\sum_{n\leq N}\bfu(n)\tend{N}{\infty}0$. In a similar way (by considering finite rotations), one can deduce $\frac1N\sum_{n\leq N}\bfu(an+b)\tend{N}{\infty}0$. Thus,~\eqref{eq:defMOMOSI} can be seen as a form of aperiodicity. A further analysis reveals that, in fact, we deal with a special behaviour of $\bfu$ on a typical short interval. All strongly aperiodic multiplicative functions satisfy~\eqref{eq:Mobius-like} (this follows from Theorem~A.1~\cite{MR3435814}), hence condition~\eqref{eq:Mobius-like} is satisfied both for $\mob$ and~$\lio$, cf.\ Section~\ref{szort}.

Recently, in~\cite{Gomilko:ab}, the strong $\bfu$-OMO property was rephrased in the language of functional analysis; and it is equivalent to
$$
\lim_{K\to \infty}\frac{1}{b_{K+1}}\sum_{k\leq K}\left\| \sum_{b_k\leq n<b_{k+1}}\bfu(n)f\circ T^n\right\|_{C(X)}=0 \text{ for all }f\in C(X),\;(b_k)\text{ as above}.
$$

Usefulness of the strong MOMO concept is seen in the following result:
\begin{Prop}[\cite{MR3874857}]\label{p:sMOMO}
If $(R,Z,\mathcal{D},\kappa)$ is an ergodic (measure-theoretic) dynamical system and $(T,X)$ is its uniquely ergodic model satisfying the strong MOMO property then {\bf all} uniquely ergodic models of $(R,Z,\mathcal{D},\kappa)$ are M\"obius orthogonal. In fact, the strong MOMO holds in all of them.
\end{Prop}

\subsection{M\"obius orthogonality of positive entropy systems}
If we take a positive entropy system $(T,X)$, it is natural to expect that it is not M\"obius orthogonal. Indeed, trivially, the full shift on $\{0,1\}^{\Z}$ is not, and more generally subshifts of finite type are not, see~\cite{Kar2}. One can also show that the subshift $(S,X_{\mob^2})$  (which is of positive entropy, see~\cite{MR3430278}) is not M\"obius orthogonal, despite the fact that $\mob^2$ itself is a completely deterministic point and M\"obius orthogonality holds at it: $\lim_{N\to\infty}\frac1N\sum_{n\leq N}f(S^n\mob^2)\mob(n)=0$ for each $f\in C(X_{\mob^2})$, see~\cite{MR3821717}.
However, Sarnak's conjecture does not exclude a possibility that {\bf some} positive entropy system is also M\"obius orthogonal.\footnote{It is mentioned in~\cite{Sa} that Bourgain (unpublished) had such a construction.}
Downarowicz and Serafin proved the following general result:
\begin{Th}[\cite{MR3961705}]\label{dose1}
Fix an integer $N\geq2$. Let $\bfu$ be any bounded, real, aperiodic sequence. Then, there
exists a subshift $(S,X)$ over $N$ symbols of entropy arbitrarily close to $\log N$, uncorrelated
to $\bfu$: $\lim_{N\to\infty}\frac1N\sum_{n\leq N}f(S^nx)\bfu(n)=0$ for each $f\in C(X)$ and $x\in X$.\end{Th}

Even more surprisingly, they proved a uniform version of the above result:
\begin{Th}[\cite{Do-Se1902.04162}] \label{dose2}
Under the same assumption on $\bfu$, given $N\geq2$, there exists a strictly ergodic subshift over $N$ symbols, of entropy arbitrarily close to $\log N$,  {\bf uniformly} uncorrelated to $\bfu$.
\end{Th}

Realizing that, one might be anxious what finally is the class of systems which are M\"obius orthogonal, and in particular, why zero entropy should play a special role. As we will see however, positive entropy systems are not expected to enjoy the strong MOMO property, cf.\  Theorem~\ref{th:3wki}. Indeed, the following has been proved in~\cite{MR3874857}:\footnote{This result has been proved in \cite{MR3874857} for the Liouville function but it can also be proved for $\mob$.}
\begin{Th}[\cite{MR3874857}]\label{th:MOMOzero} Let $\bfu\in\{-1,0,1\}^{\Z}$ be a generic point for the measure $\widehat{\nu}_{\mob^2}$. Then the following conditions are equivalent:
\begin{enumerate}[(i)]
\item $(T,X)$ satisfies strong $\bfu$-OMO property.
\item $(T,X)$ is of zero entropy.
\end{enumerate}
\end{Th}

As an immediate consequence, we have the following:
\begin{Cor}\label{SMposent} If Chowla conjecture holds then a system $(T,X)$ has the strong MOMO property if and only if it has zero entropy.\end{Cor}

\subsection{\eqref{slog} vs.\ \eqref{clog}}
The logarithmic version of Sarnak's conjecture was formulated in~\cite{MR3569059} along with \eqref{clog} and it postulates that
\begin{equation}\tag{{\bf S$_{\log}$}}\label{slog}
\lim_{N\to\infty}\frac1{\log N}\sum_{n\leq N}\frac1n f(T^nx)\mob(n)=0
\end{equation}
(with all parameters as in~\eqref{e1}).
In~\cite{MR3676413}, Tao showed the following:
\begin{Th}\label{SClog}
\eqref{slog} is equivalent to \eqref{clog}.
\end{Th}
\begin{Remark}\label{tauw}
Combining Theorem~\ref{SClog} with Corollary~\ref{wnio}, we obtain that \eqref{slog} implies \eqref{chsfor} along a subsequence of logarithmic density 1. In particular, ({\bf S}) implies \eqref{chsfor} along a subsequence of full logarithmic density.
\end{Remark}
Let us here recall one more  ``logarithmic conjecture'' from~\cite{MR3676413} which confirms a special role played by nil-systems in dynamics. Let $(T_{g_0},G/\Gamma)$ be a nilrotation. Let $f\in C(G/\Gamma)$ be Lipschitz continuous and $x_0\in G$. Then (for $H\leq N$ with $H\to\infty$)
\begin{equation}\label{Slognil}\tag{{\bf S}$_{\log}^{\text{nil}}$}
\sum_{n\leq N}\frac{\sup_{g\in G}\left|\sum_{h\leq H}f(T_g^{h+n}(x_0\Gamma))\mob(n+h)\right|}n={\rm o}(H\log N).
\end{equation}
\begin{Th}[\cite{MR3676413}]
\eqref{Slognil} is equivalent to \eqref{slog} (and \eqref{clog}).
\end{Th}

Finally, as a consequence of the result on logarithmic correlations of $\mob$ of order~2 (using Corollary~\ref{lebesgue}), we obtain:
\begin{Cor}
All topological dynamical systems whose all invariant measures yield systems with singular spectrum are logarithmically M\"obius orthogonal.
\end{Cor}
In general, we do not know if we can replace ``all invariant measures'' with ``all ergodic invariant measures'' in the above corollary (the same applies to Corollary~\ref{AA}). This replacement is possible however, when there are only countably many ergodic invariant measures, cf.\ the discussion in~\cite{MR3779960} in Section~\ref{fusy}.

\subsection{\eqref{sar} vs.\ \eqref{slog}}
Clearly, \eqref{sar} implies \eqref{slog}.
As for the other direction, we have the following:
\begin{Th}[\cite{Gomilko:ab}]
Suppose that \eqref{slog} holds. Then there exists a sequence of logarithmic density 1, along which \eqref{sar} holds for all zero entropy topological dynamical systems.
\end{Th}
The idea of the proof is to use Theorem~\ref{SClog}, Remark~\ref{tauw} and then repeat the arguments from the proof of Theorem~\ref{ChtS}.

Notice that the sequence of logarithmic density~1 in the above result is universal for all zero entropy systems. In~\cite{Gomilko:ab}, one more version of M\"obius orthoginality is studied, namely so-called {\em logarithmic strong MOMO property} (cf.\ Section~\ref{smomo}):
$$
\lim_{K\to\infty}\frac{1}{\log b_{K+1}}\sum_{k\leq K}\left\|\sum_{b_k\leq n<b_{k+1}}\frac{\mob(n)}{n} f\circ T^n \right\|_{C(X)}=0.
$$
Equivalently, for all increasing sequences $(b_k)\subset \N$ with $b_{k+1}-b_k\to\infty$, all $(x_k)\subset X$ and $f\in C(X)$,
$$
\lim_{K\to \infty}\frac{1}{\log b_{K+1}}\sum_{k\leq K}\left|\sum_{b_k\leq n<b_{k+1}}\frac{1}{n}f(T^{n-b_k}x_k)\mob(n) \right|=0.
$$
\begin{Th}[\cite{Gomilko:ab}]
Assume that a topological system $(T,X)$ satisfies the logarithmic strong MOMO property. Then there exists a sequence $A=A(T,X)\subset \N$ with full logarithmic
density such that, for each $f\in C(X)$,
$$
\lim_{A\ni N\to\infty}\left\| \frac1N\sum_{n\leq N}\mob(n)f\circ T^n\right\|_{C(X)}=0.
$$
In particular, M\"obius orthogonality holds along  a subsequence of full logarithmic density.
\end{Th}

\begin{Remark} In \cite{Gomilko:ab}, using \cite{MR3779960}, it is proved that each system $(T,X)$ for which $M^e(T,X)$ is countable satisfies the logarithmic strong MOMO property, hence, for each such system Sarnak's conjecture holds in (logarithmic) density, cf.\ Theorem~\ref{HF}.\end{Remark}

\subsection{Strategies}\label{se:strategie}
The first years of activity around Sarnak's conjecture were devoted to proving M\"obius orthogonality in {\bf selected classes} of zero entropy dynamical systems. While this proved fruitful, and often some brilliant arguments were found ad hoc, with a strong dependence on the class under consideration, it quickly became clear that it won't be sufficient. Two main strategies to attack~\eqref{sar} arose:
\paragraph{A}
The first strategy is to look for some additional intrinsic structure in zero entropy systems that could be used to prove orthoginality from $\mob$, namely {\bf internal disjointness}. Here, a priori, one does not use any other property of $\mob$ than multiplicativity and boundedness.
\paragraph{B}
As we have seen, \eqref{sar} is intimately related to \eqref{chsfor}, and therefore one cannot expect to confirm~\eqref{sar} without using further {\bf number-theoretic properties of $\mob$}. This directs attention to {\bf aperiodicity} and behaviour on so-called {\bf short intervals}. It extends further to studying {\bf Furstenberg systems of $\mob$} (including the logarithmic ones) and trying to interpret arithmetic properties of $\mob$ as ergodic properties of the corresponding dynamical systems. One can hope finally to deduce (some kind of) Furstenberg (!)  disjointness of Furstenberg systems of $\mob$ with a wide subclass of zero entropy systems (hopefully, with all such systems).

 As we will see, these two approaches often intertwine, proving once again that number theory and ergodic theory should not be studied separately from each other.

\section{Arithmetic properties of the M\"obius function}\label{se:arithm}

\subsection{Multiplicativity}
\paragraph{Internal disjointness}
Joinings (introduced in a seminal paper of Furstenberg~\cite{Fu}) have been present in ergodic theory for over 50 years. Disjointness (absence of non-trivial joinings), as a form of an extremal non-isomorphism and a measure-theoretic invariant, has always played a crucial role in classification problems.\footnote{Recall also that different powers for a typical automorphism of a standard Borel space are pairwise disjoint~\cite{Ju}. See also more recent~\cite{Kanigowski_2020}.}  It appeared however in many other contexts, including homogenous dynamics, with applications in number theory. Sarnak's conjecture gave yet a new impetus, in particular for studying (approximate) disjointness for different sub-actions.

A basic method to prove orthogonality with a multiplicative function comes from the Multiplicative Orthogonality Criterion (MOC):
\begin{Th}[\cite{Ka,Bo-Sa-Zi}]\label{t:kbsz}
Assume that $(f_n)\subset\C$ is a bounded sequence.
Assume that for all (sufficiently large) prime numbers $p\neq q$,
\begin{equation}\label{kbsz1}
\lim_{N\to\infty}\frac1N\sum_{n\leq N}f_{pn}\overline{f}_{qn}=0.
\end{equation}
Then, for each bounded multiplicative function $\bfu$, we have
$
\lim_{N\to\infty}\frac1N\sum_{n\leq N}f_n\bfu(n)=0$. In particular, $(f_n)$ is M\"obius orthogonal.
\end{Th}
\begin{Remark}
Notice that Theorem~\ref{t:kbsz} does not require from $\bfu$ anything but multiplicativity and boundedness.
\end{Remark}
In the dynamical context $(T,X)$ the simplest way to use Theorem~\ref{t:kbsz} is to take $f_n=f(T^nx)$. In this form, MOC appeared for the first time in~\cite{Bo-Sa-Zi} and was used to prove that the horocycle flows are M\"obius orthogonal.
To see how MOC is used and how it is related to Furstenberg's disjointness theory~\cite{Fu}, assume 
 that $M(T,X)=\{\mu\}$, $\int_Xf\,d\mu=0$, and the corresponding measure-theoretic system is totally ergodic. Then, any measure $\rho\in V((x,x))$ (considered in the topological dynamical system $(T^p\times T^q,X\times X)$) is a joining of $T^p$ and $T^q$. If we now assume that $(T^p,X,\mu)$ and $(T^q,X,\mu)$ are disjoint for sufficiently large primes $p\neq q$ then $\rho=\mu\ot\mu$ and, as a result, the limit in~\eqref{kbsz1} equals $\int_{X\times X}f\ot\ov{f}\,d\rho=0$, i.e.\ the assumptions of MOC are satisfied.

In general, a use of MOC is not that simple. Consider an irrational rotation $Tx=x+\alpha$ on the circle $X=\R/\Z$. To see that~\eqref{kbsz1} holds for all characters, one uses the Weyl criterion on uniform distribution. However, there are continuous zero mean functions for which~\eqref{kbsz1} fails~\cite{Ku-Le}, which shows clearly, that in general we can only expect~\eqref{kbsz1} to hold for a linearly dense set of continuous functions.

In some cases, MOC cannot be applied directly (e.g.\ when the systems under consideration fail to be weakly mixing) and the spectral approach can help. Examples can be found in~\cite{Ab-Le-Ru,Bo1,Ab-Ka-Le}.

\paragraph{AOP property}
The following ergodic counterpart of MOC was developed in~\cite{Ab-Le-Ru2}: an ergodic automorphism $(T,X,\mathscr{B},\mu)$ is said to have {\em asymptotically orthogonal powers} (AOP)  if for each $f,g\in L^2_0(X,\mathscr{B},\mu)$, we have
\begin{equation}\label{momoe4}
\lim_{\PP\ni p,q\to\infty, p\neq q} \sup_{\kappa\in J^e(T^p,T^q)}\left|\int_{X\times X} f\ot g\,d\kappa\right|=0.
\end{equation}
Clearly, if the powers of $T$ are pairwise disjoint then $T$ enjoys the AOP property. However, this condition is not necessary, the powers of $T$ having AOP property may even be isomorphic. Moreover, AOP implies total ergodicity and zero entropy~\cite{Ab-Le-Ru2}. The relation between strong MOMO and AOP properties is described by the following result:
\begin{Th}[\cite{Ab-Le-Ru2,Ab-Ku-Le-Ru}]\label{thmB}
Let $\bfu$ be a bounded multiplicative function. Suppose that $(R,Z,\mathscr{C},\kappa)$ satisfies AOP. Then the following are equivalent:
\begin{itemize}
\item
$\bfu$ satisfies \eqref{eq:Mobius-like};
\item
The strong $\bfu$-OMO property is satisfied in each uniquely ergodic model $(T,X)$ of $R$.
\end{itemize}
In particular, if the above holds, for each $f\in C(X)$, we have
$$
\frac1N\sum_{n\leq N}f(T^nx)\bfu(n)\tend{N}{\infty} 0 \text{ uniformly on } X.
$$
\end{Th}

\subsection{Aperiodicity}
As all periodic sequences are orthogonal to $\mob$, one can expect that sequences with some properties similar to periodicity will also be M\"obius orthogonal. Notice also that M\"obius orthogonality of periodic sequences (\eqref{sar} for rotations on finite groups) corresponds to $\mob$ being aperiodic. This is the simplest situation where some additional properties of $\mob$ (other than multiplicativity) begin to play a significant role.
\paragraph{Zero entropy continuous interval maps} In~\cite{Kar1}, \eqref{sar} for zero entropy continuous intervals maps and orientation-preserving circle homeomorphisms is established. The starting point for developing the main tools is the result of Davenport~\eqref{daven}, which shows clearly that the examples under consideration are indeed ``relatives'' of irrational rotations. Additionally, in order to treat the case of interval maps one studies $\omega$-limit sets and it turns out that, in fact, one deals with an odometer.

\paragraph{Synchronized automata}
In~\cite{De-Dr-Mu}, Deshouillers, Drmota and M\"ullner prove that~\eqref{sar} is true for automatic sequences generated by synchronizing automata (the inputs are read with the most significant digit first). In fact, they prove orthogonality of such sequences from any bounded  function $\bfu$ that is aperiodic.

\paragraph{Almost periodic sequences}
We say that a sequence is {\em Weyl rationally almost periodic} (WRAP) whenever it can be approximated arbitrarily well by periodic sequences in Weyl pseudo-metric $d_W$ given by $d_W(x,y)=\limsup_{N\to\infty}\sup_{\ell\geq 1}\frac1N|\{ \ell \leq n\leq\ell+N : x(n)\neq y(n)\}|$. It is proved in~\cite{MR3989121} that each subshift $(S,X_x)$ given by a Weyl almost periodic sequence is M\"obius orthogonal (in fact, we have orthogonality to any bounded aperiodic arithmetic function $\bfu$).

\subsection{Behaviour on short intervals}\label{szort}
During the last four years, an enormous progress concerning the short interval behavior of strongly aperiodic multiplicative functions  has been made due to the breakthrough result of Matom\"aki and Radziwi\l\l. Their main result of~\cite{Ma-Ra}, for $\mob$, in its simplified form can be written as
\begin{equation}\label{condmu}
\lim_{\substack{M,H\to \infty\\ \text{ with }H=o(M)}}\frac1M\sum_{1\leq m\leq M}\frac1H\left| \sum_{m\leq h<m+H}\mob(h)\right|=0.
\end{equation}
This gave an impetus to study convergence on short intervals in ergodic theory and it has become a new, crucial player from the point of view of Sarnak's conjecture. Condition~\eqref{condmu} can be also reformulated in the following way: for each $(b_n)\subset\N$ with $b_{n+1}-b_n\to\infty$, we have
\begin{equation}\label{12}
\lim_{K\to\infty}\frac{1}{b_{K+1}}\sum_{k\leq K}\left|\sum_{b_k\leq n<b_{k+1}}\mob(n) \right|=0,
\end{equation}
cf.\ Section~\ref{smomo}.

\paragraph{Almost periodic sequences}
In~\cite{MR3989121}, in case of WRAP $x$, the authors also ask about the behaviour of averages of the form
\begin{equation}\label{eq:sc-ep-7-a}
\frac1H\sum_{m\leq h<m+H}f(S^hz)\mob(n)
\end{equation}
(where $z\in X_x$) for large values of $H$ and arbitrary $m\in\N$. Under \eqref{chsfor}, convergence to zero uniformly in $m$ does not take place, however, it is shown in~\cite{MR3989121} that for a ``typical'' $m\in\N$ the averages in \eqref{eq:sc-ep-7-a} are small. The key argument in the proof comes from a result of Matom\"aki, Radziwi\l \l \ and Tao:
\begin{Th}[\cite{MR3435814}]
For each periodic sequence $a(n)$, we have
\begin{equation}\label{sarSHORT}
\lim_{\substack{H,M\to\infty\\ H={\rm o}(M)}}\frac1M\sum_{M\leq m<2M}\left|\frac1H\sum_{m\leq h<m+H}a(h)\boldsymbol{\mu}(h)\right|=0.
\end{equation}
\end{Th}
As a consequence, we have:
\begin{Th}[\cite{MR3989121}]\label{Wshort-2}
Suppose that $x\in \mathbb{A}^\Z$ is WRAP. Then for all $f\in C(X_x)$ and $z\in X_x$,
\begin{equation}
\label{eq:sc-ep-4}
\lim_{\substack{H,M\to\infty\\ H={\rm o}(M)}} \frac1M\sum_{M\leq m<2M}\Big|\frac1H\sum_{m\leq h<m+H}f(S^hz)\boldsymbol{\mu}(h)\Big|=0.
\end{equation}
\end{Th}
Moreover, it is shown that all synchronized automata yield WRAP sequences. Thus, the above theorem strengthens the aforementioned  result  by Deshouillers, Drmota and M\"ullner in~\cite{De-Dr-Mu}.

\paragraph{Rigid systems}
In~\cite{Kanigowski:aa}, Kanigowski, Lema\'{n}czyk and Radziwi\l{}\l{} study rigid systems.\footnote{A measure-theoretic system $(R,Z,\mathscr{C},\kappa)$ is {\em rigid} if, for some increasing sequence $(q_n)$ of natural numbers, we have $f\circ R^{q_n}\to f$ in $L^2(Z,\kappa)$ for each $f\in L^2(Z,\kappa)$. Rigid systems are of zero entropy. Moreover, the typical measure-theoretic automorphism is rigid and weakly mixing.} To formulate their results, we need some definitions and facts. Given a natural number $q$, the sum $\sum_{\PP\ni p | q} 1/p$ is called the {\em prime volume} of $q$. The prime volume grows slowly with $q$:
$$
\sum_{p|q}\frac1p\leq \log\log\log q+O(1).
$$
However ``most'' of the time, the prime volume of $q$ stays bounded: if we set
$$
\mathcal{D}_j:=\Big\{q\in \N : \sum_{p|q}\frac{1}{p}<j\Big\},
$$
then
$d(\mathcal{D}_j)\to1$ when $j\to\infty$.
A topological system $(T,X)$  is said to be {\em good} if for every $\nu\in M(X,T)$ at least one of the following conditions holds:
\begin{itemize}
\item
({\bf BPV rigidity}): $(T, X, \mathscr{B}, \nu)$ is rigid along a sequence $(q_n)_{n\geq1}$ with {\bf bounded prime volume}, i.e.\ there exists $j$ such that $(q_n)_{n\geq1}\subset \mathcal{D}_j$;
\item
({\bf PR rigidity}): $(T,X,\mathscr{B},\nu)$ has {\bf polynomial rate} of rigidity, i.e.\ there exists a  linearly dense (in $C(X)$) set $\mathcal{F}\subset C(X)$ such that for each $f\in\mathcal{F}$ we can find $\delta>0$ and a sequence $(q_n)_{n\geq1}$ satisfying
$$
\sum_{j=-q_n^\delta}^{q_n^\delta}\|f\circ T^{jq_n}-f\|_{L^2(\nu)}^2\to 0.
$$
\end{itemize}
They prove the following:
\begin{Th}
\begin{enumerate}[(a)]
\item
Assume that $(T,X)$ is a topological system such that $(T,X,\mathscr{B}(X),\mu)$ is good. Then $(T,X)$ is M\"obius orthogonal.
\item
Suppose that each {\bf ergodic} invariant measure of $(T,X)$ yields either BPV rigidity or PR rigidity and $M^e(T,X)$ is countable then $(T,X)$ is M\"obius orthogonal.
\end{enumerate}
\end{Th}
A key tool here is a strengthening of the main result of Matom\"aki and Radziwi\l{}\l{}~\cite{Ma-Ra} (cf.\ \eqref{condmu}) to short interval behaviour along arithmetic progressions:
\begin{Th}[\cite{Kanigowski:aa}] \label{klr}
For each $\varepsilon>0$, there exists $L_0$ such that for each $L\geq L_0$ and $q\geq 1$ satisfying $\sum_{p|q}1/p\leq(1-\varepsilon)\sum_{p\leq L}1/p$ we can find $M_0=M_0(q,L)$ such that for all $M\geq M_0$, we have
$$
\sum_{j=0}^{M/Lq}\sum_{a=0}^{q-1}\left|\sum_{\substack{m\in[z+jLq,z+(j+1)Lq] \\ m\equiv a\bmod q}}\mob(m) \right|<\varepsilon M
$$
for some $0\leq z<Lq$.
\end{Th}

\begin{Remark} Despite the fact that PR rigidity does not seem to be stable under different (uniquely ergodic) models of a measure-preserving transformation, assuming $(T,X)$ is uniquely ergodic, it is proved that $(T,X)$ satisfies the  strong MOMO property whenever its unique invariant measure yields either BPV or PR rigidity. Via Proposition~\ref{p:sMOMO}, we obtain that if in a model PR rigidity holds, all of the models are M\"obius orthogonal. This, in particular, applies to all ergodic transformations with discrete spectrum.
Moreover, it is shown in \cite{Kanigowski:aa}  that for a.e.\ IET (of $d\geq 3$ intervals) BPV rigidity holds, so a.e.\ IET (and all their uniquely ergodic models) is M\"obius orthogonal. This is to be compared with previously known results for 3-IETs~\cite{Bo1,Chaika_2019,Ferenczi_2018,Karagulyan:aa}.
Other applications are given for $C^{2+\vep}$ Anzai skew products and for some so-called Rokhlin extensions of rotations.
\end{Remark}

One more of the consequences is the following result:
\begin{Cor}[\cite{Kanigowski:aa}]\label{norigid}
No Furstenberg system of the M\"obius function $\mob$ is either BPV or PR rigid. The same holds for the Liouville function $\lio$.
\end{Cor}

\subsection{Logarithmic Furstenberg systems}\label{fusy}
Frantzikinakis and Host study logarithmic Furstenberg systems associated to $\mob$ (and $\lio$). They prove the following remarkable result:
\begin{Th}[\cite{MR3779960}]\label{HF}
Each zero entropy topological system $(T,X)$ with only countably many ergodic measures is logarithmically M\"obius orthogonal.
\end{Th}
In particular, uniquely ergodic systems of zero topological entropy satisfy~\eqref{slog}. The key argument in the proof of Theorem~\ref{HF} is the following structural result on the logarithmic Furstenberg systems of $\mob$ and $\lio$:
\begin{Th}[\cite{MR3779960}]\label{HF1}
Each logarithmic Furstenberg system of $\mob$ or $\lio$ is a factor of a system that:
\begin{itemize}
\item has no irrational spectrum,
\item has ergodic components isomorphic to direct products of infinite-step nilsystems and Bernoulli systems.
\end{itemize}
\end{Th}
The starting point for the proof of the above theorem, resulting in a reduction of the problem to purely ergodic context, is an identity of Tao (implicit in~\cite{MR3569059}) showing that self-correlations of $\mob$ (and $\lio$) are  averages of its dilated self-correlations with prime dilates. Frantzikinakis and Host also prove that logarithmic Furstenberg systems of $\mob$ (and $\lio$) are ``almost determined'' by strongly stationary processes (introduced by Furstenberg and Katznelson in the 90's). The structure of (measure-theoretic) dynamical systems given by strongly stationary processes has been described by Jenvey~\cite{Jenvey_1997} who proved that ergodicity implies Bernoulli (cf.\  Theorem~\ref{t:fr1} and~\ref{t:fr1a}) and Frantzikinakis~\cite{Frantzikinakis_2004} who described the ergodic decomposition in the non-ergodic case.

The above results are extended in~\cite{Frantzikinakis:aa} to strongly aperiodic multiplicative functions. Moreover, the following multi-dimensional result is proved:
\begin{Th}\label{nooo}
Let $f_1,\dots, f_\ell\colon\N\to\U$ be multiplicative functions. Let $(R,Y)$ be a topological dynamical system and let $y\in Y$ be a logarithmically generic point for a measure $\nu$ with zero entropy and having at most countably many ergodic components, all of which are totally ergodic. Then for every $g\in C(Y)$ that is orthogonal in $L^2(\nu)$ to all $R$-invariant functions we have
\begin{equation}
\lim_{N\to\infty}\frac{1}{\log N}\sum_{n\leq N}\frac{g(R^ny)\prod_{j=1}^{\ell}f_j(n+h_j)}{n}=0
\end{equation}
for all $h_1,\dots, h_\ell\in\Z$.
\end{Th}
The unweighted version of (with $g=1$) is expected to hold if the shifts are distinct and at least one of the multiplicative functions is strongly aperiodic. This is the logarithmically averaged version of Elliott conjecture~\cite{MR1042765,MR1222182,MR3435814}. In the special case of irrational rotations, with $\ell=1$, Theorem~\ref{nooo} is the logarithmically averaged variant of a classical result of Daboussi~\cite{MR675168,MR0332702,AST_1975__24-25__321_0}. Already in case $\ell=2$ it is completely new.

\section{Future Directions}\label{se:future}
\subsection{Detecting zero entropy}
As shown in Theorem~\ref{th:MOMOzero} and Corollary~\ref{SMposent}, under the Chowla conjecture, the M\"obius function is a sequence that:
\begin{itemize}
\item is strong MOMO orthogonal to all zero entropy systems,
\item is never strong MOMO orthogonal to any positive entropy system.
\end{itemize}

In view of this, it is natural to study the following problem:
\begin{Question}
Which numerical sequences distinguish between zero and positive entropy systems?
\end{Question}
Note that in view of the results of Downarowicz and Serafin (Theorem~\ref{dose1} and Theorem~\ref{dose2}), ``usual orthogonality'' or even its uniform version is insufficient for these needs.

\subsection{Proving the strong MOMO property}
It was already asked in \cite{MR3821717} whether whenever we have a zero entropy system $(T,X)$ for which we can prove M\"obius orthogonality, then we can prove the strong MOMO property. Recently, Lema\'{n}czyk and M\"ullner in~\cite{Lema_czyk_2020} proved the strong MOMO property for (primitive) automatic sequences (previously known to be M\"obius orthogonal by~\cite{M_llner_2017}), answering a question from~\cite{MR3821717}. Yet another question persists:
\begin{Question} Do horocycle flows satisfy the strong MOMO property?\end{Question}
We do not even know  whether M\"obius orthogonality takes place in all uniquely ergodic models of horocycle flows.

\subsection{Mixing properties of Furstenberg systems} Corollary~\ref{norigid} induces the following problem:
\begin{Question}
Are Furstenberg systems of $\lio$ mildly mixing?
\end{Question}
For $\mob$ we need to take into account that its Furstenberg systems have the discrete spectrum factor given by the Mirsky measure of $\mob^2$.

\subsection{Furstenberg disjointness in non-ergodic case}
As Host and Frantzikinakis' analysis shows, if the (potential) logarithmic Furstenberg systems of $\lio$ or $\mob$ are non-ergodic, then they are very non-ergodic. One of open questions by Frantzikinakis is whether the system $\mathbb{T}^2\ni (x,y)\mapsto (x,y+x)\in \mathbb{T}^2$ considered with Lebesgue measure can be a Furstenberg system of $\lio$. Of course this example is a measure-theoretic system which is Furstenberg disjoint from all ergodic systems. It seems to be a problem of independent interest to fully understand the class of transformations disjoint from all ergodic transformations.

\section*{Acknowledgements}
\noindent
Research supported by Narodowe Centrum Nauki grant UMO-2019/33/B/ST1/00364. 

\nociteNew{Ab-Le-Ru1}
\nociteNew{Baake:2015aa}
\nociteNew{MR3731019}
\nociteNew{MR3803141}
\nociteNew{MR1954690}
\nociteNew{Bo0}
\nociteNew{MR3296562}
\nociteNew{MR3296562}
\nociteNew{Da-Te}
\nociteNew{Do-Ka}
\nociteNew{Dr}
\nociteNew{MR3825824}
\nociteNew{Fe-Ku-Le-Ma}
\nociteNew{MR3850672}
\nociteNew{Gr}
\nociteNew{MR3821719}
\nociteNew{MR3947636}
\nociteNew{Keller:2017aa}
\nociteNew{MR3803667}
\nociteNew{Ku-Le-We}
\nociteNew{Ku-Le-We1}
\nociteNew{MR4000514}
\nociteNew{Li-Sa}
\nociteNew{MATOM_KI_2016}
\nociteNew{Ma-Ri2}
\nociteNew{MR3612882}
\nociteNew{MR0021566}
\nociteNew{Mi}
\nociteNew{Pe-Hu}
\nociteNew{Ry}
\nociteNew{MR3666035}
\nociteNew{MR3660308}
\nociteNew{MR3829173}
\nociteNew{McNamara:aa}
\nociteNew{HLSY}
\nociteNew{MR3859364}
\nociteNew{Sa:Af}
\nociteNew{Sawin:aa}
\nociteNew{Ei1}
\nociteNew{Tablog5}
\nociteNew{Veech_2017}
\nociteNew{MR3820018}
\nociteNew{Sun:aa}
\nociteNew{MR3819999}
\nociteNew{Forni:aa}
\nociteNew{Konieczny_2020}
\nociteNew{MR3927855}
\nociteNew{Houcein-el-Abdalaoui:aa}
\nociteNew{He:aa}
\nociteNew{El_Abdalaoui_2018}

\allowdisplaybreaks
\small

\bibliographystyle{abbrv}
\bibliography{sarnak-ency}

\begin{thebibliography}{10}

\bibitem{Houcein-el-Abdalaoui:aa}
E.~H. {\noopsort{Abdalaoui}}El~Abdalaoui.
\newblock Oscillating sequences, {G}owers norms and {S}arnak's conjecture.
\newblock Preprint, \url{https://arxiv.org/abs/1704.07243}.

\bibitem{El_Abdalaoui_2018}
E.~H. {\noopsort{Abdalaoui}}El~Abdalaoui, G.~Askri, and H.~Marzougui.
\newblock M{\"o}bius disjointness conjecture for local dendrite maps.
\newblock {\em Nonlinearity}, 32(1):285--300, Dec 2018.

\bibitem{Ab-Le-Ru1}
E.~H. {\noopsort{Abdalaoui}}El~Abdalaoui, M.~Lema\'nczyk, and T.~de~la Rue.
\newblock A dynamical point of view on the set of {$\mathscr{B}$}-free
  integers.
\newblock {\em Int. Math. Res. Not.}, 2015(16):7258--7286, 2015.

\bibitem{Baake:2015aa}
M.~Baake and C.~Huck.
\newblock Ergodic properties of visible lattice points.
\newblock {\em Proceedings of the Steklov Institute of Mathematics},
  288(1):165--188, 2015.

\bibitem{MR1954690}
V.~Bergelson and I.~Ruzsa.
\newblock Squarefree numbers, {IP} sets and ergodic theory.
\newblock In {\em Paul {E}rd{\H o}s and his mathematics, {I} ({B}udapest,
  1999)}, volume~11 of {\em Bolyai Soc. Math. Stud.}, pages 147--160. J{\'a}nos
  Bolyai Math. Soc., Budapest, 2002.

\bibitem{Bo0}
J.~Bourgain.
\newblock M{\"o}bius-{W}alsh correlation bounds and an estimate of {M}auduit
  and {R}ivat.
\newblock {\em J. Anal. Math.}, 119:147--163, 2013.

\bibitem{MR3927855}
R.~{\noopsort{Breteche}}de~la Bret\`eche and G.~Tenenbaum.
\newblock A remark on {S}arnak's conjecture.
\newblock {\em Q. J. Math.}, 70(1):371--378, 2019.

\bibitem{MR3296562}
F.~Cellarosi and I.~Vinogradov.
\newblock Ergodic properties of {$k$}-free integers in number fields.
\newblock {\em J. Mod. Dyn.}, 7(3):461--488, 2013.

\bibitem{Da-Te}
C.~Dartyge and G.~Tenenbaum.
\newblock Sommes des chiffres de multiples d'entiers.
\newblock {\em Ann. Inst. Fourier (Grenoble)}, 55(7):2423--2474, 2005.

\bibitem{Do-Ka}
T.~Downarowicz and S.~Kasjan.
\newblock Odometers and {T}oeplitz systems revisited in the context of
  {S}arnak's conjecture.
\newblock {\em Studia Math.}, 229(1):45--72, 2015.

\bibitem{Dr}
M.~Drmota.
\newblock Subsequences of automatic sequences and uniform distribution.
\newblock In {\em Uniform distribution and quasi-{M}onte {C}arlo methods},
  volume~15 of {\em Radon Ser. Comput. Appl. Math.}, pages 87--104. De Gruyter,
  Berlin, 2014.

\bibitem{MR3825824}
M.~Drmota, C.~M\"{u}llner, and L.~Spiegelhofer.
\newblock M\"{o}bius orthogonality for the {Z}eckendorf sum-of-digits function.
\newblock {\em Proc. Amer. Math. Soc.}, 146(9):3679--3691, 2018.

\bibitem{MR3731019}
A.~Dymek.
\newblock Automorphisms of {T}oeplitz {$\mathscr{B}$}-free systems.
\newblock {\em Bull. Pol. Acad. Sci. Math.}, 65(2):139--152, 2017.

\bibitem{MR3803141}
A.~Dymek, S.~Kasjan, J.~Ku\l{}aga-Przymus, and M.~Lema\'{n}czyk.
\newblock $\mathscr{B}$-free sets and dynamics.
\newblock {\em Trans. Amer. Math. Soc.}, 370(8):5425--5489, 2018.

\bibitem{Ei1}
T.~Eisner.
\newblock A polynomial version of {S}arnak's conjecture.
\newblock {\em C. R. Math. Acad. Sci. Paris}, 353(7):569--572, 2015.

\bibitem{MR3819999}
A.-H. Fan and Y.~Jiang.
\newblock Oscillating sequences, {MMA} and {MMLS} flows and {S}arnak's
  conjecture.
\newblock {\em Ergodic Theory Dynam. Systems}, 38(5):1709--1744, 2018.

\bibitem{Fe-Ku-Le-Ma}
S.~Ferenczi, J.~Ku{\l}aga-Przymus, M.~Lema{\'n}czyk, and C.~Mauduit.
\newblock Substitutions and {M}{\"o}bius disjointness.
\newblock In {\em Ergodic Theory, Dynamical Systems, and the Continuing
  Influence of John C. Oxtoby}, volume 678 of {\em Contemp. Math.}, pages
  151--173. Amer. Math. Soc., Providence, RI, 2016.

\bibitem{MR3850672}
L.~Flaminio, K.~Fr\k{a}czek, J.~Ku\l{}aga-Przymus, and M.~Lema\'{n}czyk.
\newblock Approximate orthogonality of powers for ergodic affine unipotent
  diffeomorphisms on nilmanifolds.
\newblock {\em Studia Math.}, 244(1):43--97, 2019.

\bibitem{Forni:aa}
G.~Forni and A.~Kanigowski.
\newblock Mutliple mixing and disjointness for time changes of bounded-type
  {H}eisenberg nilflows.
\newblock {\em ournal de l'{\'E}cole polytechnique --- Math{\'e}matiques},
  7:63--91, 2020.

\bibitem{MR3829173}
N.~Frantzikinakis.
\newblock An averaged {C}howla and {E}lliott conjecture along independent
  polynomials.
\newblock {\em Int. Math. Res. Not.}, 2018(12):3721--3743, 2018.

\bibitem{Gr}
B.~Green.
\newblock On (not) computing the {M}{\"o}bius function using bounded depth
  circuits.
\newblock {\em Combin. Probab. Comput.}, 21(6):942--951, 2012.

\bibitem{He:aa}
X.~He and Z.~Wang.
\newblock M{\"o}bius disjointness for nilsequences along short intervals.
\newblock Preprint, \url{https://arxiv.org/abs/1905.02864}.

\bibitem{HLSY}
W.~Huang, Z.~Lian, S.~Shao, and X.~Ye.
\newblock Sequences from zero entropy noncommutative toral automorphisms and
  {S}arnak {C}onjecture.
\newblock {\em J. Differential Equations}, 263(1):779--810, 2017.

\bibitem{MR3821719}
C.~Huck.
\newblock On the logarithmic probability that a random integral ideal is
  {$\mathcal{A}$}-free.
\newblock In S.~Ferenczi, J.~Ku\l{}aga-Przymus, and M.~Lema\'{n}czyk, editors,
  {\em Ergodic theory and dynamical systems in their interactions with
  arithmetics and combinatorics}, volume 2213 of {\em Lecture Notes in Math.},
  pages 249--258. Springer, Cham, 2018.

\bibitem{MR3947636}
S.~Kasjan, G.~Keller, and M.~Lema\'{n}czyk.
\newblock Dynamics of {$\mathscr{B}$}-free sets: a view through the window.
\newblock {\em Int. Math. Res. Not.}, 2019(9):2690--2734, 2019.

\bibitem{Keller:2017aa}
G.~Keller.
\newblock Generalized heredity in $\mathscr{B}$-free systems.
\newblock Preprint, \url{https://arxiv.org/abs/1704.04079}.

\bibitem{Konieczny_2020}
J.~Konieczny.
\newblock M{\"o}bius orthogonality for q-semimultiplicative sequences.
\newblock {\em Monatshefte f{\"u}r Mathematik}, 192(4):853--882, Jun 2020.

\bibitem{MR3803667}
J.~Konieczny, M.~Kupsa, and D.~Kwietniak.
\newblock Arcwise connectedness of the set of ergodic measures of hereditary
  shifts.
\newblock {\em Proc. Amer. Math. Soc.}, 146(8):3425--3438, 2018.

\bibitem{MR4000514}
J.~Ku\l{}aga-Przymus and M.~Lema\'{n}czyk.
\newblock M\"{o}bius disjointness along ergodic sequences for uniquely ergodic
  actions.
\newblock {\em Ergodic Theory Dynam. Systems}, 39(10):2793--2826, 2019.

\bibitem{Ku-Le-We}
J.~Ku\l{}aga-Przymus, M.~Lema{\'n}czyk, and B.~Weiss.
\newblock On invariant measures for {$\mathscr{B}$}-free systems.
\newblock {\em Proc. Lond. Math. Soc. (3)}, 110(6):1435--1474, 2015.

\bibitem{Ku-Le-We1}
J.~Ku\l{}aga-Przymus, M.~Lema{\'n}czyk, and B.~Weiss.
\newblock Hereditary subshifts whose simplex of invariant measures is
  {P}oulsen.
\newblock In {\em Ergodic theory, dynamical systems, and the continuing
  influence of {J}ohn {C}. {O}xtoby}, volume 678 of {\em Contemp. Math.}, pages
  245--253. Amer. Math. Soc., Providence, RI, 2016.

\bibitem{Li-Sa}
J.~Liu and P.~Sarnak.
\newblock The {M}{\"o}bius function and distal flows.
\newblock {\em Duke Math. J.}, 164(7):1353--1399, 2015.

\bibitem{MATOM_KI_2016}
K.~Matom\"{a}ki, M.~Radziwi\l{}\l{}, and T.~Tao.
\newblock Sign patterns of the {L}iouville and {M}\"obius functions.
\newblock {\em Forum of Mathematics, Sigma}, 4, 2016.

\bibitem{Ma-Ri2}
C.~Mauduit and J.~Rivat.
\newblock Prime numbers along {R}udin--{S}hapiro sequences.
\newblock {\em J. Eur. Math. Soc. (JEMS)}, 17(10):2595--2642, 2015.

\bibitem{McNamara:aa}
R.~McNamara.
\newblock Sarnak's conjecture for sequences of almost quadratic word growth.
\newblock Preprint, \url{https://arxiv.org/abs/1901.06460}.

\bibitem{MR3612882}
M.~K. Mentzen.
\newblock Automorphisms of subshifts defined by {$\mathcal{B}$}-free sets of
  integers.
\newblock {\em Colloq. Math.}, 147(1):87--94, 2017.

\bibitem{MR0021566}
L.~Mirsky.
\newblock Note on an asymptotic formula connected with {$r$}-free integers.
\newblock {\em Quart. J. Math., Oxford Ser.}, 18:178--182, 1947.

\bibitem{Mi}
L.~Mirsky.
\newblock Arithmetical pattern problems relating to divisibility by {$r$}th
  powers.
\newblock {\em Proc. London Math. Soc. (2)}, 50:497--508, 1949.

\bibitem{MR3859364}
R.~Peckner.
\newblock M\"{o}bius disjointness for homogeneous dynamics.
\newblock {\em Duke Math. J.}, 167(14):2745--2792, 2018.

\bibitem{Pe-Hu}
P.~A.~B. Pleasants and C.~Huck.
\newblock Entropy and diffraction of the {$k$}-free points in {$n$}-dimensional
  lattices.
\newblock {\em Discrete Comput. Geom.}, 50(1):39--68, 2013.

\bibitem{Ry}
V.~V. Ryzhikov.
\newblock Bounded ergodic constructions, disjointness, and weak limits of
  powers.
\newblock {\em Trans. Moscow Math. Soc.}, pages 165--171, 2013.

\bibitem{Sa:Af}
P.~Sarnak.
\newblock Mobius randomness and dynamics.
\newblock {\em Not. S. Afr. Math. Soc.}, 43(2):89--97, 2012.

\bibitem{Sawin:aa}
W.~Sawin.
\newblock Dynamical models for {L}iouville and obstructions to further progress
  on sign patterns.
\newblock preprint, \url{https://arxiv.org/abs/1809.03280}.

\bibitem{MR3666035}
K.~Soundararajan.
\newblock The {L}iouville function in short intervals.
\newblock {\em Ast\'{e}risque}, 2015/2016(390):Exp. No. 1119, 453--479, 2017.
\newblock S\'{e}minaire Bourbaki. Vol. 2015/2016. Expos\'{e}s 1104--1119.

\bibitem{Sun:aa}
W.~Sun.
\newblock Sarnak's conjecture for nilsequences on arbitrary number fields and
  applications.
\newblock Preprint, \url{https://arxiv.org/abs/1902.09712}.

\bibitem{Tablog5}
T.~Tao.
\newblock Furstenberg limits of the {L}iouville function.
\newblock Webpage,
  \url{https://terrytao.wordpress.com/2017/03/05/furstenberg-limits-of-the-liouville-function/}.

\bibitem{MR3820018}
J.~Ter\"{a}v\"{a}inen.
\newblock On binary correlations of multiplicative functions.
\newblock {\em Forum Math. Sigma}, 6:e10, 41, 2018.

\bibitem{Veech_2017}
W.~A. Veech.
\newblock M{\"o}bius orthogonality for generalized {M}orse-{K}akutani flows.
\newblock {\em American Journal of Mathematics}, 139(5):1157--1203, 2017.

\bibitem{MR3660308}
Z.~Wang.
\newblock M\"{o}bius disjointness for analytic skew products.
\newblock {\em Invent. Math.}, 209(1):175--196, 2017.

\end{thebibliography}


\begin{thebibliography}{10}

\bibitem{Ab-Ka-Le}
E.~H. {\noopsort{Abdalaoui}}El~Abdalaoui, S.~Kasjan, and M.~Lema{\'n}czyk.
\newblock 0-1 sequences of the {T}hue-{M}orse type and {S}arnak's conjecture.
\newblock {\em Proc. Amer. Math. Soc.}, 144(1):161--176, 2016.

\bibitem{Ab-Ku-Le-Ru}
E.~H. {\noopsort{Abdalaoui}}El~Abdalaoui, J.~Ku\l{}aga-Przymus,
  M.~Lema{\'n}czyk, and T.~de~la Rue.
\newblock The {C}howla and the {S}arnak conjectures from ergodic theory point
  of view.
\newblock {\em Discrete Contin. Dyn. Syst.}, 37(6):2899--2944, 2017.

\bibitem{MR3874857}
E.~H. {\noopsort{Abdalaoui}}El~Abdalaoui, J.~Ku\l{}aga-Przymus,
  M.~Lema\'{n}czyk, and T.~de~la Rue.
\newblock M\"{o}bius disjointness for models of an ergodic system and beyond.
\newblock {\em Israel J. Math.}, 228(2):707--751, 2018.

\bibitem{Ab-Le-Ru}
E.~H. {\noopsort{Abdalaoui}}El~Abdalaoui, M.~Lema\'nczyk, and T.~de~la Rue.
\newblock On spectral disjointness of powers for rank-one transformations and
  {M}{\"o}bius orthogonality.
\newblock {\em J. Funct. Anal.}, 266(1):284--317, 2014.

\bibitem{Ab-Le-Ru2}
E.~H. {\noopsort{Abdalaoui}}El~Abdalaoui, M.~Lema\'nczyk, and T.~de~la Rue.
\newblock Automorphisms with {Q}uasi-discrete {S}pectrum, {M}ultiplicative
  {F}unctions and {A}verage {O}rthogonality {A}long {S}hort {I}ntervals.
\newblock {\em Int. Math. Res. Not.}, 2017(14):4350--4368, 2017.

\bibitem{MR3989121}
V.~Bergelson, J.~Ku\l{}aga-Przymus, M.~Lema\'{n}czyk, and F.~K. Richter.
\newblock Rationally almost periodic sequences, polynomial multiple recurrence
  and symbolic dynamics.
\newblock {\em Ergodic Theory Dynam. Systems}, 39(9):2332--2383, 2019.

\bibitem{Bo1}
J.~Bourgain.
\newblock On the correlation of the {M}oebius function with rank-one systems.
\newblock {\em J. Anal. Math.}, 120:105--130, 2013.

\bibitem{Bo-Sa-Zi}
J.~Bourgain, P.~Sarnak, and T.~Ziegler.
\newblock Disjointness of {M}\"obius from horocycle flows.
\newblock In {\em From {F}ourier analysis and number theory to {R}adon
  transforms and geometry}, volume~28 of {\em Dev. Math.}, pages 67--83.
  Springer, New York, 2013.

\bibitem{Chaika_2019}
J.~Chaika and A.~Eskin.
\newblock M{\"o}bius disjointness for interval exchange transformations on
  three intervals.
\newblock {\em Journal of Modern Dynamics}, 14(1):55--86, 2019.

\bibitem{Ch}
S.~Chowla.
\newblock {\em The {R}iemann hypothesis and {H}ilbert's tenth problem}.
\newblock Mathematics and Its Applications, Vol. 4. Gordon and Breach Science
  Publishers, New York, 1965.

\bibitem{AST_1975__24-25__321_0}
Collectif.
\newblock Fonctions multiplicatives presque p\'eriodiques {B}.
\newblock In {\em Journ\'ees arithm\'etiques de Bordeaux}, number 24-25 in
  Ast\'erisque, pages 321--324. Soci\'et\'e math\'ematique de France, 1975.

\bibitem{MR0332702}
H.~Daboussi and H.~Delange.
\newblock Quelques propri\'{e}t\'{e}s des fonctions multiplicatives de module
  au plus \'{e}gal {\`a} {$1$}.
\newblock {\em C. R. Acad. Sci. Paris S\'{e}r. A}, 278:657--660, 1974.

\bibitem{MR675168}
H.~Daboussi and H.~Delange.
\newblock On multiplicative arithmetical functions whose modulus does not
  exceed one.
\newblock {\em J. London Math. Soc. (2)}, 26(2):245--264, 1982.

\bibitem{Da}
H.~Davenport.
\newblock On some infinite series involving arithmetical functions. {II}.
\newblock {\em Quart. J. Math. Oxford}, 8:313--320, 1937.

\bibitem{De-Dr-Mu}
J.-M. Deshouillers, M.~Drmota, and C.~M\"ullner.
\newblock Automatic sequences generated by synchronizing automata fulfill the
  {S}arnak conjecture.
\newblock {\em Studia Math.}, 231(1):83--95, 2015.

\bibitem{MR2809170}
T.~Downarowicz.
\newblock {\em Entropy in dynamical systems}, volume~18 of {\em New
  Mathematical Monographs}.
\newblock Cambridge University Press, Cambridge, 2011.

\bibitem{MR3961705}
T.~Downarowicz and J.~Serafin.
\newblock Almost {F}ull {E}ntropy {S}ubshifts {U}ncorrelated to the
  {M}\"{o}bius {F}unction.
\newblock {\em Int. Math. Res. Not.}, 2019(11):3459--3472, 2019.

\bibitem{Do-Se1902.04162}
T.~Downarowicz and J.~Serafin.
\newblock A strictly ergodic, positive entropy subshift uniformly uncorrelated
  to the {M}\"obius function.
\newblock {\em Studia Math.}, 2019.
\newblock Published online (Online First).

\bibitem{MR1042765}
P.~D. T.~A. Elliott.
\newblock Multiplicative functions {$|g|\leq 1$} and their convolutions: an
  overview.
\newblock In {\em S\'{e}minaire de {T}h\'{e}orie des {N}ombres, {P}aris
  1987--88}, volume~81 of {\em Progr. Math.}, pages 63--75. Birkh\"{a}user
  Boston, Boston, MA, 1990.

\bibitem{MR1292619}
P.~D. T.~A. Elliott.
\newblock On the correlation of multiplicative functions.
\newblock {\em Notas Soc. Mat. Chile}, 11(1):1--11, 1992.

\bibitem{MR1222182}
P.~D. T.~A. Elliott.
\newblock On the correlation of multiplicative and the sum of additive
  arithmetic functions.
\newblock {\em Mem. Amer. Math. Soc.}, 112(538):viii+88, 1994.

\bibitem{Ferenczi_1997}
S.~Ferenczi.
\newblock Measure-theoretic complexity of ergodic systems.
\newblock {\em Israel Journal of Mathematics}, 100(1):189--207, Dec 1997.

\bibitem{MR3821717}
S.~Ferenczi, J.~Ku\l{}aga-Przymus, and M.~Lema\'{n}czyk.
\newblock Sarnak's conjecture: what's new.
\newblock In S.~Ferenczi, J.~Ku\l{}aga-Przymus, and M.~Lema\'{n}czyk, editors,
  {\em Ergodic theory and dynamical systems in their interactions with
  arithmetics and combinatorics}, volume 2213 of {\em Lecture Notes in Math.},
  pages 163--235. Springer, Cham, 2018.

\bibitem{Ferenczi_2018}
S.~Ferenczi and C.~Mauduit.
\newblock On {S}arnak's conjecture and {V}eech's question for interval
  exchanges.
\newblock {\em Journal d'Analyse Math{\'e}matique}, 134(2):545--573, Feb 2018.

\bibitem{Frantzikinakis_2004}
N.~Frantzikinakis.
\newblock The structure of strongly stationary systems.
\newblock {\em Journal d'Analyse Math{\'e}matique}, 93(1):359--388, Dec 2004.

\bibitem{MR3742396}
N.~Frantzikinakis.
\newblock Ergodicity of the {L}iouville system implies the {C}howla conjecture.
\newblock {\em Discrete Anal.}, pages Paper No. 19, 41, 2017.

\bibitem{MR3779960}
N.~Frantzikinakis and B.~Host.
\newblock The logarithmic {S}arnak conjecture for ergodic weights.
\newblock {\em Ann. of Math. (2)}, 187(3):869--931, 2018.

\bibitem{Frantzikinakis:aa}
N.~Frantzikinakis and B.~Host.
\newblock Furstenberg systems of bounded multiplicative functions and
  applications.
\newblock To appear in Int. Math. Res. Not.,
  \url{https://arxiv.org/abs/1804.08556}.

\bibitem{Fu}
H.~Furstenberg.
\newblock Disjointness in ergodic theory, minimal sets, and a problem in
  {D}iophantine approximation.
\newblock {\em Math. Systems Theory}, 1:1--49, 1967.

\bibitem{MR3821718}
A.~Gomilko, D.~Kwietniak, and M.~Lema\'{n}czyk.
\newblock Sarnak's conjecture implies the {C}howla conjecture along a
  subsequence.
\newblock In S.~Ferenczi, J.~Ku\l{}aga-Przymus, and M.~Lema\'{n}czyk, editors,
  {\em Ergodic theory and dynamical systems in their interactions with
  arithmetics and combinatorics}, volume 2213 of {\em Lecture Notes in Math.},
  pages 237--247. Springer, Cham, 2018.

\bibitem{Gomilko:ab}
A.~Gomilko, M.~Lema\'{n}czyk, and T.~de~la Rue.
\newblock M\"obius orthogonality in density for zero entropy dynamical systems.
\newblock To appear in Pure and Applied Functional Analysis,
  \url{https://arxiv.org/abs/1905.06563}.

\bibitem{Gr-Ta}
B.~Green and T.~Tao.
\newblock The {M}{\"o}bius function is strongly orthogonal to nilsequences.
\newblock {\em Ann. of Math. (2)}, 175(2):541--566, 2012.

\bibitem{Huang:aa}
W.~Huang, J.~Liu, and K.~Wang.
\newblock M{\"o}bius disjointness for skew products on a circle and a
  nilmanifold.
\newblock Preprint, \url{https://arxiv.org/abs/1907.01735}.

\bibitem{HUANG2019827}
W.~Huang, Z.~Wang, and X.~Ye.
\newblock Measure complexity and {M}\"obius disjointness.
\newblock {\em Advances in Mathematics}, 347:827 -- 858, 2019.

\bibitem{MR3959363}
W.~Huang, Z.~Wang, and G.~Zhang.
\newblock M\"{o}bius disjointness for topological models of ergodic systems
  with discrete spectrum.
\newblock {\em J. Mod. Dyn.}, 14:277--290, 2019.

\bibitem{Jenvey_1997}
E.~Jenvey.
\newblock Strong stationarity and de {F}inetti's theorem.
\newblock {\em Journal d'Analyse Math{\'e}matique}, 73(1):1--18, Dec 1997.

\bibitem{Ju}
A.~{\noopsort{Junco}}del~Junco.
\newblock Disjointness of measure-preserving transformations, minimal
  self-joinings and category.
\newblock In {\em Ergodic theory and dynamical systems, {I} ({C}ollege {P}ark,
  {M}d., 1979--80)}, volume~10 of {\em Progr. Math.}, pages 81--89.
  Birkh{\"a}user, Boston, Mass., 1981.

\bibitem{Kam}
T.~Kamae.
\newblock Subsequences of normal sequences.
\newblock {\em Israel J. Math.}, 16:121--149, 1973.

\bibitem{Kanigowski:aa}
A.~Kanigowski, M.~Lema{\'n}czyk, and M.~Radziwi{\l}{\l}.
\newblock Rigidity in dynamics and {M}\"obius disjointness.
\newblock Preprint, \url{https://arxiv.org/abs/1905.13256}.

\bibitem{Kanigowski_2020}
A.~Kanigowski, M.~Lema{\'n}czyk, and C.~Ulcigrai.
\newblock On disjointness properties of some parabolic flows.
\newblock {\em Inventiones mathematicae}, 221(1):1--111, Jan 2020.

\bibitem{Kar1}
D.~Karagulyan.
\newblock On {M}{\"o}bius orthogonality for interval maps of zero entropy and
  orientation-preserving circle homeomorphisms.
\newblock {\em Ark. Mat.}, 53(2):317--327, 2015.

\bibitem{Kar2}
D.~Karagulyan.
\newblock On {M}{\"o}bius orthogonality for subshifts of finite type with
  positive topological entropy.
\newblock {\em Studia Math.}, 237(3):277--282, 2017.

\bibitem{Karagulyan:aa}
D.~Karagulyan.
\newblock Hausdorff dimension of a class of three-interval exchange maps.
\newblock {\em Discrete and Continuous Dynamical Systems}, 40(3):1257--1281,
  2020.

\bibitem{Ka}
I.~{\noopsort{Katai}}K{\'a}tai.
\newblock A remark on a theorem of {H}. {D}aboussi.
\newblock {\em Acta Math. Hungar.}, 47(1-2):223--225, 1986.

\bibitem{Ku-Le}
J.~Ku\l{}aga-Przymus and M.~Lema{\'n}czyk.
\newblock The {M}{\"o}bius function and continuous extensions of rotations.
\newblock {\em Monatsh. Math.}, 178(4):553--582, 2015.

\bibitem{Lema_czyk_2020}
M.~Lema{\'n}czyk and C.~M{\"u}llner.
\newblock Automatic sequences are orthogonal to aperiodic multiplicative
  functions.
\newblock {\em Discrete Contin. Dyn. Syst.}, 40(12):6877--6918, 2020.

\bibitem{Ma-Ra}
K.~Matom\"aki and M.~Radziwi{\l}{\l}.
\newblock Multiplicative functions in short intervals.
\newblock {\em Ann. of Math. (2)}, 183(3):1015--1056, 2016.

\bibitem{MR3435814}
K.~Matom\"{a}ki, M.~Radziwi\l{}\l{}, and T.~Tao.
\newblock An averaged form of {C}howla's conjecture.
\newblock {\em Algebra Number Theory}, 9(9):2167--2196, 2015.

\bibitem{M_llner_2017}
C.~M{\"u}llner.
\newblock Automatic sequences fulfill the {S}arnak conjecture.
\newblock {\em Duke Mathematical Journal}, 166(17):3219--3290, Nov 2017.

\bibitem{MR3430278}
R.~Peckner.
\newblock Uniqueness of the measure of maximal entropy for the squarefree flow.
\newblock {\em Israel J. Math.}, 210(1):335--357, 2015.

\bibitem{MR3810678}
M.~Ram~Murty and A.~Vatwani.
\newblock A remark on a conjecture of {C}howla.
\newblock {\em J. Ramanujan Math. Soc.}, 33(2):111--123, 2018.

\bibitem{MR3821722}
O.~Ramar\'{e}.
\newblock Chowla's conjecture: from the {L}iouville function to the {M}oebius
  function.
\newblock In S.~Ferenczi, J.~Ku\l{}aga-Przymus, and M.~Lema\'{n}czyk, editors,
  {\em Ergodic theory and dynamical systems in their interactions with
  arithmetics and combinatorics}, volume 2213 of {\em Lecture Notes in Math.},
  pages 317--323. Springer, Cham, 2018.

\bibitem{Sa}
P.~Sarnak.
\newblock Three lectures on the {M}{\"o}bius function, randomness and dynamics.
\newblock Webpage, \url{http://publications.ias.edu/sarnak/}.

\bibitem{MR3569059}
T.~Tao.
\newblock The logarithmically averaged {C}howla and {E}lliott conjectures for
  two-point correlations.
\newblock {\em Forum Math. Pi}, 4:e8, 36, 2016.

\bibitem{MR3676413}
T.~Tao.
\newblock Equivalence of the logarithmically averaged {C}howla and {S}arnak
  conjectures.
\newblock In {\em Number theory---{D}iophantine problems, uniform distribution
  and applications}, pages 391--421. Springer, Cham, 2017.

\bibitem{Ta5}
T.~Tao.
\newblock The logarithmically averaged and non-logarithmically averaged
  {C}howla conjectures.
\newblock Webpage,
  \url{https://terrytao.wordpress.com/2017/10/20/the-logarithmically-averaged-and-non-logarithmically-averaged-chowla-conjectures/}.

\bibitem{MR3938639}
T.~Tao and J.~Ter\"{a}v\"{a}inen.
\newblock Odd order cases of the logarithmically averaged {C}howla conjecture.
\newblock {\em J. Th\'{e}or. Nombres Bordeaux}, 30(3):997--1015, 2018.

\bibitem{MR4039498}
T.~Tao and J.~Ter\"{a}v\"{a}inen.
\newblock The structure of correlations of multiplicative functions at almost
  all scales, with applications to the {C}howla and {E}lliott conjectures.
\newblock {\em Algebra Number Theory}, 13(9):2103--2150, 2019.

\bibitem{MR3992031}
T.~Tao and J.~Ter\"{a}v\"{a}inen.
\newblock The structure of logarithmically averaged correlations of
  multiplicative functions, with applications to the {C}howla and {E}lliott
  conjectures.
\newblock {\em Duke Math. J.}, 168(11):1977--2027, 2019.

\bibitem{We9}
B.~Weiss.
\newblock Normal sequences as collectives.
\newblock In {\em Proc. {S}ymp. on {T}opological {D}ynamics and ergodic
  theory}. Univ. of {K}entucky, 1971.

\end{thebibliography}

\bibliographystyleNew{abbrv}
\bibliographyNew{sarnak-ency}

\end{document}